\documentclass[11 pt, reqno]{article}

\usepackage{amsthm}
\usepackage{amsmath}
\usepackage{amssymb}
\usepackage{amsfonts}
\usepackage{latexsym}
\usepackage{bbm}
\usepackage{comment}
\usepackage[noadjust]{cite}

\usepackage{chngcntr}
\counterwithin*{equation}{section}

\usepackage[letterpaper]{geometry}
\usepackage{color}
\usepackage{latexsym}

\newtheorem{ccounter}{ccounter}[section]
\newtheorem{thm}[ccounter]{Theorem}
\newtheorem{lem}[ccounter]{Lemma}
\newtheorem{cor}[ccounter]{Corollary}
\newtheorem{defn}[ccounter]{Definition}
\newtheorem{prop}[ccounter]{Proposition}
\newtheorem{ass}[ccounter]{Assumption}
\newtheorem{ex}[ccounter]{Example}

\def\bet{\begin{thm}}
\def\eet{\end{thm}}
\def\bel{\begin{lem}}
\def\eel{\end{lem}}
\def\bas{\begin{ass}}
\def\eas{\end{ass}}
\def\bec{\begin{cor}}
\def\eec{\end{cor}}
\def\bed{\begin{defn}}
\def\eed{\end{defn}}
\def\bep{\begin{prop}}
\def\eep{\end{prop}}
\def\beq{\begin{equation}}
\def\eeq{\end{equation}}
\def\bea{\begin{equation*}}
\def\eea{\end{equation*}}
\def\tr{\mathrm{tr}}
\def\bex{\begin{ex}}
\def\eex{\end{ex}}
\def\bp{\begin{proof}}
\def\ep{\end{proof} \\}

\def\1{{\mathbbm 1}}
\def\pp{\mathbb{P}}
\def\fc{\mathrm{fc}}

\def\benr{\begin{enumerate}[label=(\roman*)]}
\def\eenr{\end{enumerate}}

\def\blam{\boldsymbol{\lambda}}
\def\btlam{\tilde{\boldsymbol{\lambda}}}

\def\N{\mathbb{N}}
\def\T{\mathcal{T}}

\def\R{\mathbb{R}}

\def\e{\mathrm{e}}
\def\P{\mathbb{P}}
\def\E{\mathbb{E}}

\def\msc{m_{\mathrm{sc} } }

\def\one{{\mathbbm 1}}
\def\eps{\varepsilon}

\def\scrho{\rho_{\mathrm{sc} }}
\def\Im{\operatorname{Im}}
\def\Re{\operatorname{Re}}
\newcommand{\bma}{\begin{bmatrix}}
\newcommand{\ema}{\end{bmatrix}}

\def\tr{\operatorname{tr}}
\def\ni{\noindent}

\def\scg{\gamma}
\def\i{i}
\def\A{\mathcal{A}}
\renewcommand\d{d}
\newcommand{\llbracket}{[[}
\newcommand{\rrbracket}{]]}

\def\rr{\mathbb{R}}

\def\del{\partial}

\def\O{{O}}

\def\ee{\mathbb{E}}

\def\tillam{\tilde{\lambda}}
\def\var{\operatorname{Var}}

\date{\today}

\makeatletter
\def\blfootnote{\xdef\@thefnmark{}\@footnotetext}
\makeatother

\begin{document}



\begin{table}
\centering
\begin{tabular}{c}
\multicolumn{1}{c}{\Large{\bf Comparison theorem for some extremal eigenvalue statistics}}\\
\\
\\
\end{tabular}

\begin{tabular}{c c c c c}
Benjamin Landon$^1$ & & Patrick Lopatto$^2$ & & Jake Marcinek$^2$\\
\\
& &{\small{$^1$Department of Mathematics}}& & \\
& &{\small{Massachusetts Institute of Technology}}& & \\
 \\
& &{\small{$^2$Department of Mathematics}}& & \\
& &{\small{Harvard University}}& & \\
 \\
\\
& & \today & &\\
\\
\end{tabular}
\\
\begin{tabular}{p{14 cm}}
\small{\bf Abstract.}  We introduce a method for the comparison of some extremal eigenvalue statistics of random matrices.  For example, it allows one to compare the maximal eigenvalue gap in the bulk of two generalized Wigner ensembles, provided that the first four moments of their matrix entries match.  As an application, we extend results of Ben Arous--Bourgade and Feng--Wei that identify the limit of the maximal eigenvalue gap in the bulk of the GUE to all complex Hermitian generalized Wigner matrices. 
\end{tabular}
\end{table}


\section{Introduction}  {\blfootnote{P.L. is partially supported by the NSF Graduate 
Research Fellowship Program under Grant DGE-1144152.}There has been significant progress in understanding local eigenvalue statistics of random matrices over the past decade.  A fundamental example is the proof of  bulk universality for Wigner matrices established in the series of works \cite{EKYY12,EKYY13,EPR10,ESY11,EY15,rigidity}.  It states that in the high dimensional limit, the local eigenvalue fluctuations are universal and depend only on the symmetry class of the random matrix ensemble. Parallel results were established in certain cases in \cite{TV10,TV11}, including a ``four moment theorem," which shows that the local statistics are determined by the first four moments of the entries. By local eigenvalue statistics, we refer to the fluctuations of the eigenvalues in the bulk  (the interior of the spectrum) on the inter-particle scale $N^{-1}$.

 The \emph{generalized Wigner ensembles} constitute a fundamental class of random matrix ensembles. These are self-adjoint $N \times N$ random matrices whose entries above the diagonal are independent centered random variables.  The variances of the entries are taken to be order $N^{-1}$, with the constraint that the variances of each row sum to $1$.   
 Examples include the Gaussian Orthogonal Ensemble and Gaussian Unitary Ensemble (GOE/GUE),  whose entries are real and complex Gaussians, respectively. 
In these cases, arguments based on orthogonal polynomials show that the asymptotic local statistics are described by explicit formulas. Universality implies that these formulas describe the asymptotic local statistics all generalized Wigner ensembles.

Examples of local statistics typically studied are the distribution of a single gap between consecutive eigenvalues, the local correlation functions at an energy $E \in (-2, 2)$ in the bulk, and averaged versions of these quantities.
 In this work we consider certain extremal eigenvalue statistics. In particular, we study the maximal gap between consecutive eigenvalues in the bulk. To place our main technical contribution in context, we briefly review the steps for proving bulk universality for generalized Wigner matrices.  We refer to \cite{EYbook,survey} for a comprehensive review and bibliography.

The first step is dynamical and originates in the work of Erd\H{o}s, Schlein, and Yau \cite{ESY11}.  It is based around analyzing the local time to equilibrium of the matrix eigenvalues under a stochastic matrix dynamics known as \emph{Dyson Brownian motion} (DBM).  The equilibrium measure of this dynamics is the GOE or GUE, depending on the symmetry class under consideration. It is now known that local observables reach their equilibrium distribution under the DBM flow after only a short time.  At the level of matrices, this result proves that the local statistics of matrix ensembles with small Gaussian perturbations coincide with those of the Gaussian ensembles. It therefore establishes that bulk universality holds for the former matrices, which are known as \emph{Gaussian divisible ensembles}.

The second step, extending universality from Gaussian divisible ensembles to all generalized Wigner matrices, is a density argument. One approach is the four moment method of Tao and Vu \cite{TV11,TV10}. Their work implies that the local eigenvalue statistics of two matrix ensembles with four matching moments coincide as $N \to \infty$, which yields universality for ensembles whose moments match those of the GOE or GUE.  In general, it is possible to construct Gaussian divisible ensembles matching four moments of a given generalized Wigner ensemble. Combining this fact and the four moment method with the first step, we obtain universality for all generalized Wigner matrices. Other approaches to this step include the Green function comparison theorem \cite{EYY12}, the matrix continuity estimate of \cite{QUE}, and the reverse heat flow of \cite{EPR10}.

We  refer to the two steps of this approach as the dynamical step and the comparison step, and now consider how they may be applied to prove universality of the maximal eigenvalue gap. For the dynamical part, the works \cite{landon2019fixed,wigfixed} employ a coupling of Dyson Brownian motion to an auxiliary process which is already at equilibrium.  They show that the $i$th eigenvalue gap of the DBM is the same as the $i$th gap of the equilibrium process, down to the scale $N^{-1-c}$.  The work \cite{landon2019fixed}  proved further that this estimate holds with a probability large enough to permit a union bound over all $N$ eigenvalues.  This is sufficient to obtain universality of the maximal gap distribution for Gaussian divisible ensembles.  We remark that the polynomial error rate $c>0$ obtained in this work is not optimal, and in particular, it is not sufficient to treat the minimal eigenvalue gap.

For the comparison step, the reverse heat flow provides strong estimates, but requires that the matrix entry distributions are smooth, an assumption we would like to avoid.  The other approaches are based on the Lindeberg moment matching strategy and apply either to individual eigenvalues or resolvent entries. They are local in that they work only near a spectral energy $E$ or for functions of a few eigenvalues. However, the maximal gap is a statistic involving a macroscopic number of eigenvalues distributed throughout the spectrum. Therefore, a new approach is needed in order to handle the maximal gap and other extremal spectral statistics for ensembles with discrete distributions.

Our contribution is to provide comparison theorems for the maximal gap between consecutive eigenvalues in the bulk (and certain other extremal statistics).  These theorems extend universality for the largest eigenvalue gaps from the Gaussian divisible ensembles to general matrices.  As a consequence, we show that the results of Ben Arous--Bourgade \cite{maxgapGUE} and Feng--Wei \cite{feng2018large} on the distribution of the maximal gaps of the GUE in fact hold for all Hermitian generalized Wigner matrices; we describe these more below.  Our results can also be viewed as providing a four-moment theorem for extremal spectral statistics in the spirit of Tao and Vu \cite{TV10,TV11}.  

At a technical level, we construct a regularization of the maximal eigenvalue gap that is amenable to comparison.  The four-moment approach is based around the Lindeberg moment matching strategy, which requires estimates on the derivatives of the regularization with respect to the matrix entries.  The application of the Lindeberg strategy to random matrices originates in the work of Chatterjee  \cite{chatterjee}.  Our derivative estimates also allow us to extend the matrix continuity approach of \cite{QUE} to maximal eigenvalue gaps. 

The regularization we use for the maximum is
\beq \label{eqn:freereg}
\max_i x_i \sim \frac{1}{\beta} \log \left( \sum_i \e^{ \beta x_i} \right)
\eeq
for appropriately chosen $\beta \to \infty$ as $N\to \infty$. This regularization was used previously in random matrix theory in conjunction with the Lindeberg strategy by Korada and Montanari \cite{montlind} to show universality of the minima of loss functions associated to certain statistical algorithms; this is otherwise unrelated to our work. Our inspiration for this regularization comes from statistical mechanics, where it represents the ``zero temperature limit'' of a Gibbs measure. (In our notation, $\beta$ is the inverse temperature.) This interpretation plays a key role in many areas of mathematics inspired by statistical mechanics, for example in the study of random constraint satisfaction problems through the theory of mean field spin glasses.

We do not apply the regularization \eqref{eqn:freereg} directly to the eigenvalue gaps.  Derivatives of eigenvalues with respect to matrix entries are singular as they involve in the denominator the differences of nearby eigenvalues.  Instead, we first construct a regularized version of the eigenvalues using the matrix resolvent in combination with the Helffer--Sj\"ostrand formula.  Part of this construction is essentially implicit in the work of Knowles and Yin \cite{ky}.  Our use of the Helffer--Sj\"ostrand formula simplifies many of the estimates, and we feel there is value in isolating the relevant arguments and presenting the construction as its own lemma. Moreover, we give a precise estimate on the relation between the accuracy of the regularization and the growth of its derivatives with respect to matrix entries.  Roughly, if the regularization is chosen so that the error between it and the eigenvalue is less than $N^{-1-\delta}$, then the $k$th derivative is no larger than $N^{-1+(k-1)\delta}$.  Tracking this dependence may be useful for future works probing the eigenvalue behavior below the inter-particle scale $N^{-1}$.

We now turn to related results in the literature. The problem of determining the distribution of the largest gap between eigenvalues was raised for random unitary matrices by Diaconis \cite{diaconis2003patterns}. As motivation, he gives a striking conjecture relating the distribution of the largest gaps between the unitary eigenvalues to extremal gaps between zeros of the Riemann zeta function \cite{odlyzko1987distribution}. For more on this connection, we refer to the reader to \cite[Section 1.3]{maxgapGUE}.

A heuristic approach to studying the maximal gap of random unitary matrices was given by Vinson, who obtained the correct order of magnitude \cite{vinson2011closest}. In \cite{maxgapGUE}, Ben Arous and Bourgade  proved for both unitary matrices and the GUE that the maximal gap, normalized by $\sqrt{ \log (N)}/N$, converges to an explicit constant (which depends on the ensemble). (For the GUE, they restrict to gaps in the bulk of the spectrum.)  
Feng and Wei \cite{feng2018large} found that, up to a deterministic re-centering, the maximal gap of both ensembles fluctuates on the scale $( N \sqrt{ \log(N)} )^{-1}$ around its limit and identified the limiting distributions. We determine the maximal gap down to the scale $N^{-1-c}$ for some $c>0$, and so our work extends the GUE results to all Hermitian generalized Wigner ensembles. In another direction, we note Figalli and Guionnet extended the results of \cite{maxgapGUE} to $\beta$-ensembles with $\beta=2$ \cite{figalli2016universality}.

Our paper is related to recent work of Bourgade \cite{bourgadeprep}.  He gives a new approach to the analysis of Dyson Brownian motion and obtains strong estimates enabling access not only to the maximal gap, but also the  \emph{minimal} gap, a much more singular quantity which is beyond the methods of our work.  The work \cite{bourgadeprep} relies on the reverse heat flow for the comparison step.  Our works are therefore complementary, as we focus on proving comparison theorems and do not study the dynamical side.  In particular, our comparison theorems apply to ensembles whose matrix entries are discrete random variables, which are outside the scope of the reverse heat flow technique. 

We note that other previous works have also studied the minimal eigenvalue gap. Vinson obtained the limit of the smallest eigenvalue gap for both the Circular Unitary Ensemble and GUE in his Ph.D.\ thesis \cite{vinson2011closest}; multiple smallest gaps of these ensembles were then considered in \cite{maxgapGUE}. Feng, Tian, and Wei have established the behavior of the smallest gaps of the Gaussian Orthogonal Ensemble; this work builds on earlier work of Feng and Wei on the circular $\beta$-ensembles \cite{Feng2018small, Feng2018GOE}. Other works have studied Wigner \cite{erdos2010Wegner, nguyen2017random} and sparse \cite{luh2018sparse, lopattoprep} matrices, obtaining tail bounds for the size of the gap. Showing universality for the minimal gap of Hermitian generalized Wigner matrices without imposing a smoothness assumption on the entries remains an intriguing open problem.

Another interesting problem is to compute the normalized limit of the maximal bulk eigenvalue gap for the GOE. Unfortunately, the argument of \cite{maxgapGUE} for the GUE relies on its determinantal structure, which is not present for the GOE. If the normalized limit of the maximal bulk gap could be established for the GOE, our arguments would immediately show it is the same for all real symmetric generalized Wigner matrices. One could also ask analogous questions about extremal gaps for $\beta$-ensembles with arbitrary $\beta$ and potential.

The remainder of the paper is organized as follows. In Section \ref{sec:result} we define the eigenvalue statistics we consider and state our main comparison results and universality corollaries.  In Section \ref{sec:proofs} we prove some of our comparison theorems.  Section \ref{sec:ev} contains an eigenvalue regularization lemma based on the method of Knowles and Yin \cite{ky}.  In Section \ref{sec:int} we prove the remainder of our comparison theorems.  Section \ref{sec:univ} and the appendix are devoted to deducing the universality corollaries from existing results in the literature and our comparison results.

\vspace{5 pt}

\ni {\bf Acknowledgments.} P.L. thanks Arka Adhikari, Amol Aggarwal, and Paul Bourgade for useful discussions. The authors thank the anonymous referee for valuable suggestions on an earlier version of this paper, especially for the alternative proof of Corollary \ref{cor:univ} appearing in Section~\ref{sec:univ}.

\section{Main result} \label{sec:result}

\subsection{Definition of model}

We begin by defining generalized Wigner matrices.

\bed \label{d:gwigner} A generalized Wigner matrix $H_N$ is a real symmetric or complex Hermitian $N\times N$ matrix whose upper triangular elements $\{h_{ij}\}_{i\le j}$ are independent random variables with mean zero and variances $\sigma_{ij}^2=\E(|h_{ij}|^2)$ that satisfy 

\beq \sum_{i=1}^N \sigma_{ij}^2 =1 \quad \text{ for all } j\in \llbracket1,N\rrbracket\eeq
and
\beq \frac{c}{N} \le \sigma_{ij}^2 \le \frac{C}{N} \quad  \text{ for all } i,j \in \llbracket1,N\rrbracket\eeq
for some constants $c, C >0$. 

 When $H_N$ is Hermitian, we further assume 

\beq c \var \Re h_{ij} \le \var \Im h_{ij} \le C \var \Re h_{ij}, \quad i\neq j\eeq

\ni and that $\Re h_{ij}$ and $\Im h_{ij}$ are independent.

Finally, we suppose that the normalized entries have finite moments, uniformly in $N$, $i$, and $j$, in the sense that for all $p\in \N$ there exists a constant $\mu_p$ such that
\beq \E\left| \frac{h_{ij}}{\sigma_{ij}} \right|^p \le \mu_{p}.\eeq

\ni for all $N$, $i$, and $j$.
\eed
Throughout this work, we suppress the dependence of various constants in our results on the constants in this definition. This dependence does not affect our arguments in any substantial way.

We remark that the moment condition on the entries is technical and used only to get convergence of the largest gaps in the topology of $L^p$ for every $p$ in \eqref{e:lpconclusion} below.  If one desires only convergence in some $L^p$ for fixed $p$, then the requirement can be relaxed to requiring that the matrix entries have a certain large moment.  

We also use the following notion of \emph{overwhelming probability}.
\bed We say that a set of events $\{\mathcal A(u)\}_{u \in U^{(N)}}$, where  $U^{(N)}$ is a parameter set which may depend on $N$, holds with \emph{overwhelming probability} if, for any $D>0$, there exists $N\big(D, U^{(N)}\big)$ such that for $ N \ge N\big(D, U^{(N)}\big)$,
\begin{equation}
 \inf_{u \in U^{(N)}} \P\left( \mathcal A(u) \right) \ge 1 - N^{-D}.
\end{equation}
\eed

\subsection{Comparison results}

We first present a theorem on the comparison of the largest eigenvalue gaps of generalized Wigner matrices with matching moments. 

We  denote eigenvalues of matrices $\lambda_1 \leq \cdots \leq \lambda_N$ and use $\blam$ to denote the vector of $N$ eigenvalues.   The statistics we  consider involve the maxima of eigenvalue gaps away from the spectral edges.  Taking such a maximum requires the specification of which eigenvalues we  take the maximum over.  We  do this in two different ways.  The first statistic, which we denote by $\T_{\ell, J}$, considers gaps $\lambda_{i+1} - \lambda_i$ with $i \in J$ where $J$ is a fixed index set.  The second, which we denote by $\widehat{\T}_{\ell, I}$, considers gaps with $\lambda_i \in I$ where $I$ is an interval.  We first consider $\T$ before considering $\widehat{\T}$ later.  These statistics are closely related, but dealing with $\widehat{\T}$ is more difficult than $\T$ due to the fact that the eigenvalue fluctuations  naturally cause the set of indices over which the maximum is taken to itself be random.

 For any index set $J \subseteq \llbracket1, N-1\rrbracket$ we denote by $\T_{\ell, J} ( {\blam} )$ the $\ell$th largest eigenvalue gap in the index set $J$, formally defined by
\begin{equation*}
\sup_{j \in J } \left\{ \lambda_{j+1}-\lambda_j :  \exists j_1 \neq \cdots \neq j_{l-1} \in J \text{ and } \lambda_{j_1+1} - \lambda_{j_1} \geq  \cdots \geq \lambda_{j_{l-1}+1} - \lambda_{j_{l-1}} \geq \lambda_{j+1} - \lambda_j \right\}.
\end{equation*}

The proofs of the next two theorems appear in Section \ref{sec:proofs}.

\bet \label{thm:4mm}
Let $H^{(v)}$ and $H^{(w)}$ be two generalized Wigner matrices with matrix elements denoted by $v_{ij}$ and $w_{ij}$ respectively.  Suppose that there exists a constant $ c>0$ such that
\beq \label{eqn:mmatch}
\left| \ee[ v_{ij}^a \bar{v}_{ij}^b ]- \ee [ w^a_{ij} \bar{w}_{ij}^b ]  \right| \leq N^{-2-c}
\eeq
for all nonnegative integers $a,b$ such that $a+ b \leq 4$.  Fix $\alpha \in (0,1/2)$ and choose some index set $J\subset \llbracket\alpha N, (1-\alpha) N\rrbracket$. Let $S \in C^{\infty} ( \rr )$ be a test function, and let $\nu = N / \sqrt{ \log(N)}$.  Suppose $\ell = \ell(N)$ satisfies $\ell  = N^{a_N}$ for some sequence $a_N \ge 0$. Then there exists constants $c_1=c_1(c) > 0$ and $C=C(c,\alpha )>0$ such that if $a_N \le c_1$, then\footnote{We did not work out the optimal dependence of $c_1$ on $c$, but as indicated in the proof below, it suffices to take $c_1 \le c/100$.}
\beq \label{eqn:mmresult}
\left| \E_{H^{(v)}} {S(\nu \T_{\ell,J}(\blam))} - \E_{H^{(w)}} {S(\nu \T_{\ell,J}(\blam))} \right| \leq C \left(\sup_{0\le d \le 5} \| S^{(d)}\|_\infty \right) N^{-c_1}.
\eeq
\eet
The following interpolation, called Dyson Brownian motion, is often used in random matrix theory.   We define for any real symmetric generalized Wigner matrix $X$ the matrix Ornstein-Uhlenbeck process $X_t = \{x_{ij}(t)\}_{i,j=1}^N$ given by
\beq 
\label{d:OU} dx_{ij}(t) = \frac{dB_{ij}(t)}{\sqrt{N}} - \frac{1}{2N s_{ij}} x_{ij}(t)\,dt,
 \eeq
 where the initial data are the entries $\{x_{ij}\}$ of $X$ and $s_{ij} = \E [ x_{ij}^2 ]$. We let $\blam(t)$ denote the vector of eigenvalues of $X_t$.
 
 For a Hermitian generalized Wigner matrix, we define Dyson Brownian motion by applying \eqref{d:OU} separately to its real and imaginary parts, with corresponding rescaling parameters $s^{\mathbb R}_{ij}= \E [ ( \Re x_{ij})^2 ]$ and $s^{\mathbb C}_{ij} = \E [ ( \Im x_{ij})^2 ]$ for the $dt$ terms.
\bet \label{thm:ito} Consider the SDE \eqref{d:OU} for a matrix $X_t$ with initial data given by some generalized Wigner matrix $X$.  Suppose $t\in (N^{-1+\delta}, N^{-1/2 - \delta})$ for some $\delta >0$. Fix $\alpha \in (0,1/2)$ and choose some index set $J\subset \llbracket\alpha N, (1-\alpha) N\rrbracket$. Let $S\in C^\infty(\R)$ and set $\nu = N / \sqrt{\log N}$. Suppose $\ell = \ell(N)$ satisfies $\ell  = N^{a_N}$ for some sequence $a_N \ge 0$. Then there exist constants $c=c(\delta)>0$ and $C=C(\delta, \alpha)>0$ such that if $a_N< c$, then 
\beq  \label{eqn:itoresult} \left|  \E \left[S\left( \nu \mathcal T_{\ell, J }(\blam(t) ) \right) \right]-  \E \left[S( \nu \mathcal T_{\ell, J }(\blam(0) ) )\right] \right| \le C \left(\sup_{0\le d \le 3} \| S^{(d)}\|_\infty \right) N^{-c}.\eeq 
\eet
With more work we can extend the above comparison theorems to the following related extremal statistic.  Let $I \subseteq (-2+\kappa, 2 - \kappa )$ be an interval with $\kappa >0$.  For $v_1 \leq v_2 \cdots \leq v_N$ with $\boldsymbol v \in \rr^N$, let
\beq
\widehat{\T}_{1, I} (\boldsymbol v) = \max_{ i : v_i \in I }(  v_{i+1} - v_i) ,
\eeq
and define similarly $\widehat{\T}_{\ell, I } (\boldsymbol v)$ to be the $\ell$th largest gap $v_{i+1} - v_i$ such that $v_i \in I$.  Then we have the following analog of the above theorems. Its proof is deferred until Section \ref{sec:int}.
\bet \label{thm:hatT}
First, 
suppose that $H^{(v)}$ and $H^{(w)}$ are two generalized Wigner matrices such that the first four moments of the matrix entries match in the sense of \eqref{eqn:mmatch}.  Then the estimate \eqref{eqn:mmresult} holds for some constants $C,c_1>0$ with $\widehat{\T}_{\ell, I}$ in place of $\T_{\ell, J}$ (which depend only on $c$ from \eqref{eqn:mmatch} and $\kappa$).

Second, let $X$ be a generalized Wigner matrix and $X_t$  be the process in \eqref{d:OU}.  Let $t \in (N^{-1+\delta}, N^{1/2-\delta} )$ for some $\delta >0$.  Then the estimate \eqref{eqn:itoresult} holds with $\widehat{\T}_{\ell, I}$ in place of $\T_{\ell, J}$ for some constants $C, c >0$ (which depend only on $\delta$ and $\kappa$).

In particular, for a fixed $\kappa$, these constants are uniform in the choice of $I \subseteq (-2+\kappa, 2 - \kappa )$.
\eet

We remark that the method we introduce is fairly robust and can handle a variety of extremal spectral statistics.  For example, in light of the result \cite[Theorem 1.7]{maxgapGUE} one may wish to scale each gap by the semicircle density near $\lambda_i$, in addition to the factor $\nu$.  Our proof applies without change to the maxima over quantities such as $\{ \nu \alpha_i ( \lambda_{i+1} - \lambda_i)\}_i$ for deterministic positive $\alpha_i$ satisfying $\alpha_i \leq N^{c}$ for a sufficiently small $c>0$.  Additionally, one may also consider the deviations $\{ \nu \alpha_i ( \lambda_i - \beta_i ) \}_i$ for similar $\alpha_i$ and arbitrary constants $\beta_i$ (or their absolute values).

In certain cases it may be necessary to use the four-moment approach instead of the matrix continuity estimate, if one can prove universality only for Gaussian divisible ensembles with somewhat larger Gaussian components. This means using Theorem \ref{thm:4mm} instead of Theorem \ref{thm:ito} (or the first claim of Theorem \ref{thm:hatT} instead of the second); compare the restriction $t \ll N^{-1/2}$ with the weaker \eqref{eqn:mmatch}.

\subsection{Corollaries for the universality of extremal gaps}
The following is a corollary of Theorem \ref{thm:ito} and the homogenization result \cite[Theorem 3.1]{landon2019fixed}.  

\bec \label{cor:univ}
Let $H$ be a generalized Wigner matrix and $G$ the GOE or GUE matrix of the same symmetry class.  There are constants $C,c>0$ such that
\begin{align}  \label{eqn:univgap}
\left| \ee_H S ( \nu \T_{\ell, J} ( \blam ) ) ]- \ee_G [ S ( \nu \T_{\ell, J} ( \blam )  ] \right| \leq C N^{-c} \sup_{0 \leq d \leq 3 } ||S^{(d)}||_\infty.
\end{align}
A similar estimate holds for $\widehat{\mathcal T}_{\ell, I}$.
\eec

Together with the results of Ben Arous--Bourgade \cite[Theorem 1.7]{maxgapGUE} and Feng--Wei \cite[Theorem 2]{feng2018large} it implies the following corollary. \bec  \label{cor:herm}
Suppose $H$ is a Hermitian generalized Wigner matrix. Let $I\subset (-2,2)$ be a compact interval, and set $M  = M(I) =  \inf_I \sqrt{4-x^2}$.
\begin{enumerate}
\item Let $\ell_N = N^{k_N}$ be a sequence of positive integers with $\lim_{N\rightarrow \infty} k_N = 0$. 
Then for any $p>0$,
\beq\label{e:lpconclusion}
\frac{MN}{\sqrt{32 \log N }} \widehat{\T}_{\ell_N, I}  \xrightarrow{L^p}  1.\eeq
\item Fix an index $k$. Let 
\beq
\tau^*_k = \frac{1}{4}(2 \log (N))^{\frac{1}{2}}\left( MN \widehat{\T}_{k, I} - \sqrt{32 \log(N)} \right)+ \frac{5}{8} \log ( 2 \log (N)).
\eeq
Then for any bounded interval $I_1 \subseteq \rr_+$ we have
\beq
\pp \left( \tau^*_k \in I_1 \right) \rightarrow  \int_{I_1} \frac{ \e^{ k ( c_2 - x ) }}{(k-1)!} \e^{ - \e^{c_2 -x }}  \, d x,
\eeq
where $c_2$ is an explicit constant depending on the interval.  It is the same as in \cite[Theorem 2]{feng2018large}.
\end{enumerate}
\eec
The proofs of these corollaries appear in Section \ref{sec:univ}.

\subsection{Preliminaries}

In this section we recall some standard facts about generalized Wigner matrices.    The Stieltjes transform of $H$ is defined by
\begin{align}
m_N (z) = \frac{1}{N} \sum_i \frac{1}{ \lambda_i - z}
\end{align}
and the semicircle law and corresponding Stieltjes transform are
\begin{align}
\scrho (E) = \frac{\sqrt{ (4 - E^2)_+ }}{2 \pi }, \qquad \msc (z) = \int \frac{\scrho (x)\, dx}{ x - z }.
\end{align}
The Green's function of $H$ is
\beq
G(z) = \frac{1}{H-z},
\eeq
and its matrix elements are denoted $G_{ij} (z)$.  

We also define the spectral domain
\beq
\mathcal D = \left\{ z = E + \i \eta \in \mathbb{C} :  |E| \leq 10,  \frac{N^\delta}{N} \leq \eta \le 10 \right\}.
\eeq
We have the following from \cite{scgeneral}.

\bet[{\cite[Theorem 2.3]{scgeneral}}] \label{thm:semi}
Fix $\eps >0$ and $\delta >0$.  Then, with overwhelming probability, we have
\begin{align}\label{e:localsclaw}
\sup_{z \in \mathcal D} \left| m_N (z) - \msc (z) \right| \leq  \frac{N^{\eps}}{N \eta }.
\end{align}
For the individual Green's function elements, with overwhelming probability we have
\beq
\sup_{z \in \mathcal D} \left| G_{ij} (z) - \delta_{ij} \msc (z) \right| \leq  N^{\eps} \left( \sqrt{ \frac{ \Im  \msc (z) }{N \eta} } + \frac{1}{N \eta } \right).
\eeq
\eet

The classical eigenvalue locations of the semicircle law are denoted by $\scg_i$ and defined by
\beq
\frac{i}{N} = \int_{-2}^{\scg_i} \scrho (x) \, dx.
\eeq
We next state eigenvalue rigidity and complete eigenvector delocalization estimates. The former is {\cite[Theorem 7.6]{scgeneral}}, while the latter follows by a standard argument from the Green's function estimates in Theorem~\ref{thm:semi} (see, for example, the proof of \cite[Theorem 2.10]{lectures}).
\bet[{\cite[Theorem 7.6]{scgeneral}}]  \label{thm:rig}
The following estimates hold for a generalized Wigner matrix $H$, simultaneously for all $ i \in \llbracket1, N \rrbracket$, with overwhelming probability.  For any $\eps >0$,
\beq\label{e:rigidity}
 | \lambda_i - \scg_i | \leq  \frac{ N^{\eps}}{N^{2/3} ( \min \{ i^{1/3}, (N-i+1)^{1/3} \} ) },
\eeq
and for the eigenvector $u_i$ of $H$ corresponding to $\lambda_i$, 
\beq
|| u_i ||_\infty \leq  \frac{N^{\eps}}{\sqrt{N}}.
\eeq
\eet

\section{Maximal gap over a set of indices} \label{sec:proofs}

In this section we prove the comparison results, Theorems \ref{thm:4mm} and \ref{thm:ito}.  Both are based on defining a suitable regularization of the $k$th largest gaps and proving estimates on the partial derivatives of this regularization with respect to matrix entries.  Theorem \ref{thm:4mm} is then based on the Lindeberg strategy, while Theorem \ref{thm:ito} is based on the matrix continuity estimate of \cite[Lemma A.1]{QUE}.  The real symmetric and complex Hermitian cases are nearly identical, and we give full details in only the symmetric case as it is notationally simpler.  

Throughout this section we fix an $\alpha \in ( 0, 1/2)$ and only consider index sets $J \subset  \llbracket \alpha N, (1-\alpha ) N \rrbracket$.  We  first begin by introducing our regularization of the $L^\infty$ norm.  Given a vector $\boldsymbol v\in \R^N$, define the associated largest bulk gap by
\beq
 \mathcal T_{1} = \mathcal T_{1, J}(\boldsymbol v) = \sup_{i \colon i \in J} v_{i+1} - v_i.
 \eeq
For general $\ell$, we recall that $\mathcal T_{\ell, J}(\boldsymbol v)$ is defined as the $\ell$th largest gap $v_{i+1} - v_i$  in consecutive elements of $\boldsymbol v$  with $i \in J$ (or zero if $\ell$ is greater than the number of $i$ with $i \in J$).  It is convenient to write this quantity as
\beq 
 \mathcal T_{\ell}(\boldsymbol v) = \mathcal T_{\ell, J} (\boldsymbol v)= \sup_{i_1 < \dots < i_\ell \colon {i_k} \in J} \sum_{k=1}^\ell (v_{i_{k}+1} - v_{i_{k}} ) - \sup_{i_1 < \dots < i_{\ell-1} \colon i_k  \in J} \sum_{k=1}^{\ell -1} (v_{i_{k}+1} - v_{i_{k}} ).
 \eeq
Set $\nu = N/\sqrt{\log N}$, which represents the scale of the largest gap.  For $\beta >0$ we introduce
\beq 
Z = Z_{\ell, \beta, J} (\boldsymbol v) = \sum_{i_1 < \dots < i_\ell \colon i_k \in J} \exp \left( \beta \nu \sum_{k=1}^\ell  v_{i_k+1} - v_{i_k} \right),\quad  G = G_{\ell,\beta,J}(\boldsymbol v) = \frac{1}{\beta} \log Z,
\eeq
where we use the convention $G_{0,\beta,J}(\boldsymbol v)=0$. We also set
\beq 
F_{\ell, \beta, J}(\boldsymbol v) = G_{\ell, \beta,J}(\boldsymbol v) - G_{\ell - 1, \beta, J}(\boldsymbol v).
\eeq
The quantity $F_{\ell, \beta, J}(\boldsymbol v)$ is our regularization of the $\ell$th largest gap of the vector $\boldsymbol v$.  The following lemma is elementary and its proof is omitted.
\bel\label{l:entropy} For any $v\in \R^n$ and index set $J$, we have 
\beq | \nu \mathcal T_\ell(\boldsymbol v) -  F_{\ell,\beta}(\boldsymbol v) |  < \frac{2 \ell \log N}{\beta}.\eeq
\eel



Instead of directly taking $\boldsymbol v$ to be the eigenvalues of a random matrix, we are going to work with a regularized version, which we  denote by the vector $\btlam= ( \tillam_i)$.  The construction of this regularization is somewhat involved, but we will see that $\tillam_i$ is primarily an integral of the empirical Stieltjes transform $m_N(z)$ over a certain domain in the complex plane.  
We state the following lemma which asserts the existence of $\btlam$ and gives the key estimates on its derivatives that are used in the present section; its proof is deferred until Section~\ref{sec:ev}.  

Given a matrix $H$, $0 \leq \theta \leq 1$, and indices $c, d$, we denote by $\theta^{cd} H$ the matrix with entries
\beq \label{eqn:thetaH}
\left( \theta^{cd} H \right)_{ij} = \begin{cases} \theta H_{ij} & \mbox{if } (i, j) = (c, d) \mbox{ or } (j, i) = (c, d), \\
H_{ij} & \mbox{otherwise.} \end{cases}
\eeq
This object naturally appears as the error in Taylor expansions with respect to matrix entries.
\bel \label{lem:ev}
Fix $\delta, \eps>0$.  There are smooth functions $\tillam_i (X)$ (depending on $\eps, \delta$) on the space of $N \times N$ symmetric matrices with the following properties.  Suppose that $H$ is a real symmetric generalized Wigner matrix.  With overwhelming probability, uniformly for all $i \in \llbracket \alpha N, (1-\alpha)N \rrbracket$, all integers $1 \le k \le 5$, and all choices of indices $1 \le a,b,c,d \le N$, we have
\beq \label{eqn:evreg}
|\tillam_i (H) - \lambda_i (H) | \leq C \frac{N^{\eps}}{N^{1+\delta}}, \qquad  \sup_{ 0 \leq \theta \leq 1 }| \del^k_{ab} \tillam_i (\theta^{cd} H ) | \leq C \frac{N^{\eps+(k-1)\delta}}{N},
\eeq
where $\del_{ab}=\del_{X_{ab}}$ denotes the partial derivative with respect to the matrix entry with index $(a,b)$  and $C=C(\eps,\delta)>0$ is a constant. The matrix $\theta^{cd} H$ is defined in \eqref{eqn:thetaH}. 

Further, we have
\beq\label{e:tildeas}  \sup_{ 0 \leq \theta \leq 1 }| \del^k_{ab} \tillam_i (\theta^{cd} H ) | \leq C N^{Ck}
\eeq
almost surely for all integers $1 \le k \le 5$.
 
 In the complex Hermitian case, the same estimates hold where instead $\del_{ab}^k$ is replaced by  $\del_{ \Re[X_{ab} ] }^i\del_{\Im [ X_{ab} ] }^j$ for $i+j=k$,  and $H$ is a complex Hermitian generalized matrix. 
\eel
The following is the main technical proposition of this work.  It provides a regularization of the largest eigenvalue gaps and estimates on the derivatives of this regularization with respect to the matrix entries.  
\bep \label{prop:main}
Let $X$ be a real symmetric generalized Wigner matrix, and let $J \subseteq \llbracket \alpha N, (1-\alpha ) N \rrbracket$.  Let $\beta = N^{\gamma}$ and suppose that $\ell \leq N^{\mathfrak{a}}$ for some $\gamma, \mathfrak a >0$.  Let $\tillam_i$ and $F_{\ell, \beta, J}$ be as above, and fix $\delta, \eps >0$. Then the following statements hold for all choices of indices $1\le a,b \le N$. First, with overwhelming probability,
\beq\label{e:thefirstestimate}
| \nu \mathcal{T}_{\ell, J} ( \blam) - F_{\ell, \beta, J} ( \btlam )| \leq C \nu \frac{N^{\eps}}{N^{1+\delta}} + 2 \frac{ N^{\mathfrak{a}} \log(N) }{ N^{\gamma}}.
\eeq
Secondly, with overwhelming probability for all $1\le k \le 5$,
\beq \label{eqn:TF}
\sup_\theta \left|  \del_{ij}^k F_{\ell, \beta, J} \left[ \btlam ( \theta_{ab} H) \right] \right| \leq \frac{C N^{k \delta}}{\beta}  \left( 1 + N^{k(\eps+ \gamma+\mathfrak{a}-\delta)} \right).
\eeq
Third, we have almost surely
\beq \label{e:as} \sup_\theta \left|  \del_{ij}^k F_{\ell, \beta, J} \left[ \btlam ( \theta_{ab} H) \right]\right| \leq C N^{Ck}. \eeq
Here $C=C(\delta, \eps, \gamma, \mathfrak a, \alpha)>0$ is a constant. The statements in the complex Hermitian case are adjusted as indicated in Lemma~\ref{lem:ev}.
\eep

Before proving the above proposition we derive the following elementary estimate. 
\bel
 \label{l:fpartial} The partial derivatives of $F_{\ell, \beta} (\boldsymbol v)$ with respect to the entries of the vector $v$ satisfy
\beq 
 \sum_{\underline{j}} \left|\frac{\partial^d F(\boldsymbol v)}{\partial_{j_1}\dots \partial_{j_d}}\right| \le C_d \beta^{d-1} \nu^d \ell^d,
 \eeq
where the sum runs over all multi-indices $\underline{j} = (j_1, \cdots j_d )$ with values in $J^d$, and $\del_j = \del_{v_j}$.
\eel
\begin{proof} We have 

\beq \partial_j Z(\boldsymbol v) = \beta \nu \sum_{\underline i} \eps_{\underline i}(j) \exp \left( \beta \nu \sum_{k=1}^\ell  v_{i_{k}+1} - v_{i_k} \right), \eeq

\ni where $\eps_{\underline i}(j)$ is defined as 

\begin{equation}
\eps_{\underline i}(j) = \begin{cases} 
0 &\mbox{if } j = i_k \mbox{ and } i_{k} = i_{k-1} + 1 \mbox{ for some } k \in \llbracket1,\ell\rrbracket, \\
- 1 &\mbox{if } j = i_k \mbox{ for some } k \in \llbracket1,\ell\rrbracket, \mbox{ but $i_{k} \neq i_{k-1} + 1$, }\\
1 &\mbox{if } j = i_k+1 \mbox{ for some } k \in \llbracket1,\ell\rrbracket, \mbox{ but $i_{k+ 1} \neq i_{k} + 1$, }\\
0 &\mbox{otherwise.}
\end{cases}
\end{equation}

\ni Higher derivatives are similar, yielding

\beq \frac{\partial^d Z(\boldsymbol v)}{\partial_{j_1}\dots\partial_{j_d}} =  (\beta \nu )^d\sum_{\underline i} \eps_{\underline i}(\underline j) \exp \left( \beta \nu \sum_{k=1}^\ell  v_{i_k+1} - v_{i_k} \right), \eeq

\ni with $\eps_{\underline i}(\underline j) \in \{-1,0,1\}$ satisfying $\eps_{\underline i}(\underline{j}) = 0$ if $\underline{j} \not\subset \underline{i} \cup \underline{i} + 1$ (where we abuse notation and consider these tuples as sets). Hence, for all $d$,
\begin{align}
\sum_{\underline{j}} \left|\frac{\partial^d Z(\boldsymbol v)}{\partial_{j_1}\dots \partial_{j_d}}\right| \leq ( \beta \nu)^d \sum_{\underline{i}} \exp \left( \beta \nu \sum_{k=1}^\ell v_{i_k+1} - v_{i_k } \right) \sum_{\underline{j}} | \eps_{\underline{i}} ( \underline{j} ) | \leq ( \beta \nu)^d Z ( \boldsymbol v) (2 \ell)^d.
\end{align}
By the chain rule and the definition of $F_{\ell, \beta, J} (\boldsymbol v)$ in terms of $\log(Z)$, we then find
\begin{align}
\sum_{\underline{j}} \left|\frac{\partial^d F(\boldsymbol v)}{\partial_{j_1}\dots \partial_{j_d}}\right| \leq \frac{1}{\beta} C_d \ell^d ( \beta \nu)^d
\end{align}
for some combinatorial factor $C_d$.  
 \end{proof}

We are now prepared to prove Proposition \ref{prop:main}, Theorem \ref{thm:4mm}, and Theorem \ref{thm:ito}.

\begin{proof}[Proof of Proposition \ref{prop:main}.] The estimate \eqref{e:thefirstestimate} follows from the first estimate of \eqref{eqn:evreg} and Lemma \ref{l:entropy}.   For the second estimate, fix integers $ 1 \le d \le k \le 5$ and let $\{s_i \}_{i=1}^d$ be positive integers such that $s_1+s_2 + \cdots s_d \le k$.
 We estimate the following function evaluated at the matrix $\theta_{ab} H$:
\begin{flalign}
\left| \sum_{j_1, \cdots j_d} \frac{ \partial^d F}{ \partial j_1 \cdots \partial j_d } ( \btlam) \frac{ \partial^{s_1} \tillam_{j_1}}{ \del_{ij}^{s_1}} \cdots \frac{ \partial^{s_d} \tillam_{j_d }}{ \del_{ij}^{s_d}} \right| &\leq \sum_{j_1, \cdots j_d} \left| \frac{ \partial^d F}{ \partial j_1 \cdots \partial j_d } ( \btlam) \right| N^{d \eps+k \delta - d \delta} \frac{1}{N^d} \notag\\
&\leq \frac{C_k N^{k \delta}}{\beta} \left( \frac{ N^{\eps+\mathfrak{a} + \gamma} \nu}{N^{1+\delta}} \right)^d \leq \frac{C_k N^{k \delta}}{\beta} N^{d ( \eps+\gamma + \mathfrak{a}-\delta)}.
\end{flalign}
In the final inequality we used $\nu \leq N$.  The first inequality is by Lemma \ref{lem:ev} and the second inequality is Lemma \ref{l:fpartial}.   By the chain rule, the $k$th partial derivative of $F$ with respect to the $(i,j)$ matrix element is a linear combination of such terms, and the claim follows. The third estimate is similar, using instead the last estimate of Lemma \ref{lem:ev}. \end{proof}

\begin{proof}[Proof of Theorem \ref{thm:4mm}]
 By Proposition \ref{prop:main} it is enough to prove Theorem \ref{thm:4mm} for the function $S(F(\btlam))$ instead of $\nu \mathcal T(\blam)$, assuming the parameters $\delta, \eps, \gamma, \mathfrak a$ are chosen such that $\gamma > \mathfrak a$ and $\delta > \eps$. The comparison for $S(F(\btlam))$ is a standard application of the Lindeberg four moment matching method, which we now briefly recall. In this approach, one replaces the upper triangular matrix entries of $H^{(v)}$ by those of $H^{(w)}$ one by one.  The difference \eqref{eqn:mmresult} is then a telescoping sum with $\O(N^2)$ terms, each being the difference of $S(F(\btlam))$ evaluated at two generalized Wigner matrices differing in only in the $(a,b)$ and $(a,b)$ entries.  One then Taylor expands each term around $X_{ab}=0$ to $5$th order.  The first four orders cancel to $\O(N^{-2-c'})$ thanks to the hypothesis \eqref{eqn:mmatch}, and the $5$th order remainder term  is $\O(N^{-2-c'})$ due \eqref{eqn:TF} and \eqref{e:as}.  As there are $\O (N^2)$ terms in the telescoping sum, the total error is $\O(N^{-c'})$, for some $c' >0$, and we deduce the result.  For a detailed exposition of this method, we refer to Chapter 16 of \cite{EYbook}.
  
  For clarity, and to explain the origin of the restrictions on the constants, we compute the $4$th order error term. Specifically, we fix a matrix entry $(a,b)$ and bound the error that arises when replacing $v_{ab}$ in the $(a,b)$ and $(b,a)$ entry with $w_{ab}$ (assuming the remaining entries are some collection of $v_{ij}$ and $w_{ij}$ entries). Write $F=F(\btlam)$ and $S=S(F(\btlam))$, let $S_x$ denote $S$ evaluated with the $(a,b)$ entry equal to $x$, and denote $\partial = \partial_{ab}$. Then Taylor expanding about $0$ in the $(a,b)$ entry gives that the $4$th order error term for $S_{v_{ab}} - S_0$ is 
  \beq\label{e:4expansion}
  \frac{v_{ab}^4}{4!}\left( S^{(4)} (\partial F)^4 + 6 S^{(3)} (\partial F)^2 \partial^2 F + S^{(2)} \left( 3 \left(\partial^2 F \right)^2  + 4 \left(\partial^3F\right)  \partial F \right) + S^{(1)} \partial^{(4)} F   \right),
  \eeq
  where this expression is evaluated with the $(a,b)$ entry set to zero. Taking expectation in \eqref{e:4expansion} and subtracting the analogous expression for $S_{w_{ab}} - S_0$, we see the $w^4_{ab}$ and $v^4_{ab}$ factors are independent of the terms in the parentheses, and we make take the difference of their expectations to obtain a factor that is $O(N^{-2-c})$ by \eqref{eqn:mmatch}, for some given $c$ depending on the matrix ensembles $H^{(v)}$ and $H^{(w)}$. It then suffices to show that the sum of the terms in parentheses in \eqref{e:4expansion} is $O(N^{\mathfrak b})$ for some $\mathfrak b < c$. 
  
  We consider just the $S^{(1)} \partial^{(4)} F$ term, since the others are similar.\footnote{The essential point is that the orders of the derivatives on $F$ in each term sum to $4$.} 
  Because $S^{(1)}$ is bounded by $\| S^{(1)} \|_\infty$, it suffices to bound $\partial^{(4)} F$ in expectation. 
  If we first work on the set of the overwhelming probability where \eqref{eqn:TF} holds, then it suffices to choose $\delta, \varepsilon, \gamma, \mathfrak a$ such that $\gamma = 2 \mathfrak a$, $\delta = 2\eps$, and $c/100 > \delta, \varepsilon, \gamma, \mathfrak a$, for example.\footnote{This ensures that $\gamma > \mathfrak a$  and $\delta > \eps$, as was required by the application of  Proposition \ref{prop:main} in the beginning of this proof.} On the set where this estimate does not hold, we use \eqref{e:as}.
 We may then take $c_1 = c/100$ in the statement of Theorem \ref{thm:4mm}, and this completes the proof. 
 \end{proof}
  
\begin{proof}[Proof of Theorem \ref{thm:ito}]
  Theorem \ref{thm:ito} is a consequence of Lemma A.1 of \cite{QUE} and Proposition \ref{prop:main}. 
\end{proof}

\section{Eigenvalue regularization} \label{sec:ev}

This section is devoted to the proof of Lemma \ref{lem:ev}.   We only treat the real symmetric case in detail, as the complex Hermitian case is analogous.  We remark that the choice of $\tillam_i$ is the same in both cases, as it is just a  function of the empirical Stieltjes transform $\tr (X-z)^{-1}$.  The proof of the derivative estimates is similar in both cases, with real symmetric being notationally simpler.

In the first subsection we construct the regularized eigenvalues $\tillam_i$ and at the same time prove the first estimate of \eqref{eqn:evreg}.  In the remaining subsections we estimate derivatives.
\subsection{Construction of regularized eigenvalues} \label{sec:regdef}
Let $\{\lambda_i\}_{i=1}^N$ be the eigenvalues of a generalized Wigner matrix $H$ and let $\alpha\in (0,1/2)$ be fixed.  Fix ${i \in \llbracket \alpha N , (1-\alpha) N \rrbracket}$ in the bulk.  Let $\eps_1 >0$ be a parameter.  Let $j$ and $k$ be indices such that 

\beq i - 2 N^{\eps_1} \leq j \leq i - N^{\eps_1}, \quad i + N^{\eps_1 } \leq k \leq i + 2 N^{\eps_1}, \eeq 
and set   \beq {I  = [ \gamma_j, \gamma_k ]}. \label{d:I}\eeq Let ${N(E)  = \left| \{ j : \lambda_j < E \}  \right|}$ be the eigenvalue counting function for $X$.

Using rigidity for generalized Wigner matrices, we have with  overwhelming probability that
\beq
\lambda_i - \gamma_j = \int_{\gamma_j}^{\lambda_i } \, dE = \int_{ I} \1_{ \{ \lambda_i \leq E \} } \, dE  = \int_{I} \1_{ \{ N(E) \geq i \} } \, dE.
\eeq
Let now $r: \rr \to [0, 1]$ be a smooth function such that $r (x) = 1$ for $x \geq i$ and $r (x) = 0$ for $x \leq i -1/2$.  We take $|r'| + |r''| + |r'''| \leq C$.  Then since $N(E)$ is integer valued,
\beq
\int_{I} \1_{ \{ N(E) \geq i \} } \, dE = \int_{I} r ( N (E) ) \, dE.
\eeq

Let $\delta_1> 0$ be a parameter and set
\beq
\eta_1 = N^{-1-\delta_1}.
\eeq
For each $E \in I$ we define the smoothed out eigenvalue counting function $f_E(x)$ on the scale $\eta_1$ as follows.  We define  $f_E (x) = 1$ for $x \leq E$ and $f_E (x) = 0$ for $x \geq E + \eta_1$.  We take $|f_E^{(k)} | \leq C_k \eta_1^{-k}$ for $k=1, 2, 3$.  Now, we have
\beq
| N (E) - \tr ( f_E) | \leq \left| \{ a : \lambda_a \in [E, E + \eta_1 \} \right|,
\eeq
and so by rigidity we have with overwhelming probability
\begin{align}
\left| \int_{I} r ( N (E) ) \, dE - \int_I r ( \tr ( f_E ) ) \, dE \right| & \leq C \int_I \left| N(E) - \tr ( f_E ) \right| \, dE \notag\\
 & \leq C \int_I \left| \{ a :  \lambda_a \in [ E, E + \eta_1 ] \} \right| \, dE \notag\\
 &= C \sum_a \int_I \1_{  \{ \lambda_a \in [ E , E + \eta_1 ] \} } \, dE \notag\\
 &\leq C \eta_1 | \{ a : \lambda_a \in I + [0, \eta_1 ] \} | \leq C \eta_1 N^{\eps_1} = C \frac{N^{\eps_1 - \delta_1 }}{N}. \label{e:smoothtr}
\end{align}
Let $\eta_2 = N^{-\delta_2 }/N$ for a parameter $\delta_2 \in (0,1)$, and let $\chi$ be a smooth, symmetric cutoff function with $\chi (x) = 1$ for $|x| \leq 1$ and $\chi (x) = 0$ for $|x| > 2$.  
 By the Helffer--Sj\"ostrand formula \cite{EYbook}, we have 
\begin{align}
\tr (f_E ) &= \frac{N}{  \pi } \int_{\rr^2} \left( i f_E (e) \chi' ( \sigma) - \sigma f'_E (e) \chi' ( \sigma ) \right) m_N (e + i \sigma ) \, de \, d\sigma \notag\\
&+ \frac{ i N } { \pi } \int_{ | \sigma | > \eta_2}  \int_{\R}  f''_E (e) \sigma\chi ( \sigma )   \Im [ m_N (e + \i\sigma ) ] \, de \, d\sigma \notag\\
&+ \frac{2 \i N }{  \pi} \int_0^{\eta_2} \int_{\rr}  f''_E (e) \sigma \Im [ m_N (e + i\sigma ) ] \, de \, d\sigma.
\end{align}
Since $m(z)$ is holomorphic, we have the Cauchy--Riemann equation $\partial_x \Im m ( x + iy) = - \partial_y \Re m( x + iy)$. Using this, we integrate the second term by parts twice (first in $e$ and then in $\sigma$) to obtain
\begin{align}
\tr (f_E ) &= \frac{N}{  \pi } \int_{\rr^2} \left( i f_E (e) \chi' ( \sigma) - \sigma f'_E (e) \chi' ( \sigma ) \right) m_N (e + i\sigma ) \, de \, d\sigma \notag\\
&+ \frac{ i N } { \pi } \int_{ | \sigma | > \eta_2}  \int_\R f'_E (e)  \del_{\sigma} ( \sigma \chi ( \sigma ) )   \Re [ m_N (e + i\sigma ) ] \, de \, d\sigma \notag\\
&+ \frac{2 i N }{  \pi} \int_0^{\eta_2}  \int_{ \rr} f''_E (e) \sigma \Im [ m_N (e + i\sigma ) ] \, de \, d\sigma \notag\\
&+ \frac{2 i N }{  \pi } \int_{\rr} f'_E (e) \eta_2 \Re [ m_N (e + i\eta_2 ) ]\, de.
\end{align}
 Define
\begin{align}
F_E &= \frac{N}{  \pi } \int_{\rr^2} \left( i f_E (e) \chi' ( \sigma) - \sigma f'_E (e) \chi' ( \sigma ) \right) m_N (e + i\sigma ) \, de \, d\sigma \notag\\
&+ \frac{ 2 i N } { \pi } \int_{ | \sigma | > \eta_2}  \int_\R f'_E (e)  \del_{\sigma} ( \sigma \chi ( \sigma ) )   \Re [ m_N (e + i\sigma ) ] \, de \, d\sigma. 
\end{align}
We estimate, using the definition of $f_E$,
\begin{align}
\left| \int_{I} r ( \tr f_E ) \, dE - \int_I r ( F_E) \, dE \right|  & \leq C \int_I \int_{\rr} \int_{0}^{\eta_2} N \sigma |f''_E ( e) |  \Im [ m_N (e + i\sigma ) ] \, d\sigma \, de \, dE \label{e:fefe}\\
&+ C \int_I \int_{\rr} N \eta_2 | f'_E (e) |  | m_N (e + i\eta_2 )|  \, de \, dE \notag\\
&\leq \frac{C}{ \eta_1^2} \int_{I} \int_0^{\eta_1} \int_0^{\eta_2} N \sigma \Im [ m_N (E+e + i\sigma ) ] \, d\sigma \, de \, dE \notag\\
&+C  \frac{ \eta_2}{\eta_1}  \int_{I} \int_{0}^{\eta_1} N \eta_2  | m_N (E + e + i\eta_2 ) | \, de \, dE \notag\\
&= \frac{C}{ \eta_1^2} \int_0^{\eta_2} \int_0^{ \eta_1 } \left[ \int_I N \sigma \Im [ m_N (E + e+ i\sigma ) ] \, dE \right] \, de \, d\sigma \label{e:I1} \\
&+ C \frac{\eta_2}{\eta_1} \int_{0}^{ \eta_1} \left[ \int_I N | m_N (E + e + i\eta_2 ) |  \, dE \right] \, de. \label{eqn:reg1}
\end{align}
We now estimate the inner integral of \eqref{e:I1}.  Define $E' = E + e$.  We have
\begin{align}
\int_I \sigma N \Im [ m_N (E' + i\sigma ) ] \, dE &\le  \sum_{ a : |i-a| \leq 3 N^{\eps_1} } \eta_2 \int_{I} \frac{\sigma}{ (E' - \lambda_a )^2 + \sigma^2 } \, dE \notag\\
&+ \sum_{ a : |i-a| > 3 N^{\eps_1} } \eta_2 \int_{I} \frac{\sigma}{ (E' - \lambda_a )^2 + \sigma^2 } \, dE. 
\end{align}
We estimate the first term by 
\beq
 \sum_{ a : |i-a| \leq 3 N^{\eps_1} } \eta_2 \int_{I} \frac{\sigma}{ (E' - \lambda_a )^2 + \sigma^2 } \, dE \leq  \sum_{ a : |i-a| \leq 3 N^{\eps_1} } \eta_2 \int_{\rr} \frac{\sigma}{ (E' - \lambda_a )^2 + \sigma^2 } \, dE \leq C \eta_2 N^{\eps_1}.
\eeq

For the second term we have with overwhelming probability for $a$ with $|i-a| > 3 N^{\eps_1}$ that $|E' - \lambda_a | \geq c N^{\eps_1}/N$ and so
\begin{align}
\sum_{ a : |i-a| > 3 N^{\eps_1} } \eta_2 \int_{I} \frac{\sigma}{ (E' - \lambda_a )^2 + \sigma^2 } \, dE &\leq C\eta_2  \frac{ N \sigma}{ N^{\eps_1} } | I | N \sup_{ E'} \left[ \Im  m_N (E' + iN^{\eps_1}/N) \right] \notag\\
& \leq C \sigma \eta_2 N \leq  \eta_2 N^{-\delta_2} \le \eta_2 N^{\eps_1}  .
\end{align}
In the last inequality, we used \eqref{e:localsclaw} to bound $\Im m_N$. The term $\eqref{e:I1}$ is therefore bounded by $C N^{-1 + \eps_1 + \delta_ 1 - 2\delta_2}$.

For the term \eqref{eqn:reg1} we again write

\begin{align} \int_I N|m_N(E' + i \eta_2) | \, dE  & \le   \sum_{ a : |i-a| \leq 3 N^{\eps_1} } \int  \frac{1}{|(E' - \lambda_a) + i \eta_2|} \, dE \\
& +  \sum_{ a : |i-a| > 3 N^{\eps_1} }  \int \frac{1}{|(E' - \lambda_a) + i \eta_2|} \, dE.
\end{align}

\ni The integral in the first term can be bounded directly and we obtain for that sum the bound 

\beq  C N^{\eps_1} | \log (\eta_2)| \le C N^{\eps_1 } \log (N). \eeq 

\ni For the second term we have the bound $C N^{\eps_1} \log(N)$ by rigidity \eqref{e:rigidity}, with overwhelming probability. Using these bounds, we conclude that \eqref{eqn:reg1} satisfies the bound
\beq
C \frac{\eta_2}{\eta_1} \int_{0}^{ \eta_1} \left[ \int_I N | m_N (E + e + i\eta_2 ) |  \, dE \right] \, de  \leq C N^{ - 1 + \eps_1 - \delta_2} \log ( N ) .
\eeq
Hence, recalling \eqref{e:fefe},
\beq
\left| \int_{I} r ( \tr f_E ) \, dE - \int_I r ( F_E) \, dE \right|  \leq C \frac { N^{\eps_1 + \delta_1 - 2 \delta_2 }}{N} + C N^{ - 1 + \eps_1 - \delta_2} \log ( N ).
\eeq
Using \eqref{e:smoothtr} we have therefore proven that, with overwhelming probability,
\beq
\left| ( \lambda_i - \gamma_j ) - \int_I r ( F_E) \, dE  \right| \leq C N^{\eps_1}  \frac{ \left( N^{\delta_1 - 2 \delta_2 } + N^{- \delta_1 }  + N^{-\delta_2 } \right) }{N} \log (N).
\eeq
We define $\tillam_i$ by 
\beq\label{e:evreg1}
\tillam_i = \int_I r (F_E ) d E + \gamma_j.
\eeq
We now fix parameters $\delta_1=\delta_2 = \delta$ and $\eps_1 = \eps/2$ and observe the first estimate of \eqref{eqn:evreg} holds for the $\tillam_i$ we have constructed.

\subsection{Derivative bounds}

In this subsection we obtain estimates on derivatives of the $\tillam_i$ defined in the previous subsection with respect to the matrix entries.  We first need some estimates on the Green's function.  Since in the proof of Theorem \ref{thm:4mm} we need to bound remainder terms in Taylor expansions involving the matrices $\theta_{ab} H$ as defined above, we make the following definition.

\bed
We say that rigidity and delocalization hold for a matrix $X$ with parameters $\alpha$ and $\eps$ if
\beq
\sup_{i \in \llbracket \alpha N, (1-\alpha) N \rrbracket} | \lambda_i(X) - \scg_i | \leq  \frac{N^{\eps/10}}{N},
\eeq
and for its eigenvectors $u_i$ we have
\beq
\sup_i ||u_i ||_\infty \leq \frac{N^{\eps/10}}{\sqrt{N}}.
\eeq
\eed

\subsubsection{Green function bounds}
In this section we  consider self-adjoint matrices $X$ with a Green's function we denote by $G_{ab}$.  
We have the following two \emph{a priori} estimates. The first is immediate from \cite[Lemma 10.2]{lectures} and the spectral decomposition \cite[(2.1)]{lectures}.
\bel[{\cite[Lemma 10.2]{lectures}}] \label{lem:Gbd}
Fix $\alpha\in (0,1)$ and $\eps >0$.   Let $E \in [ \scg_{\alpha N/2}, \scg_{(1-\alpha/2) N}]$ and $\eta >0$.  If rigidity and delocalization hold for $X$, then
\beq\label{e:betterlaw}
|G_{ab} (E + i\eta ) | \leq  N^{\eps} \left( \frac{1 }{ N \eta } + 1 \right).
\eeq
\eel
We also have a better bound if we are integrating in $E$.
\bel  Fix $\alpha \in (0, 1)$ and $\eps, \eps_1 \in (0, 1/2)$ with $\eps < \eps_1/2$.  Let $I \subseteq  [ \scg_{\alpha N/2}, \scg_{(1-\alpha/2) N}]$ be an interval such that $|I| \leq N^{\eps_1}/N$.   There exists a constant $C= C(\eps, \alpha)>0$ such that if rigidity and delocalization hold for $X$ with parameters $\eps$ and $\alpha$, then for any $\eta > 0$,
\beq
\int_{I} |G_{ab} (E + i\eta ) | \, dE \leq C \frac{ N^{\eps+\eps_1}}{N} \left( 1 + | \log ( \eta ) |  \right).
\eeq
\eel
\begin{proof}  By the spectral theorem and eigenvector delocalization,
\beq
\left| G_{ab} (z) \right| = \left| \sum_n \frac{ u_n (a) u_n (b) }{ \lambda_n - z } \right|
 \leq \frac{N^{\eps/5}}{N} \sum_n \frac{1}{ | \lambda_n - z | }.
\eeq
We write
\begin{align}
\int_I | G_{ab} (E + i\eta ) | \, dE &\leq \frac{ N^{\eps/5}}{N} \sum_{ n : |n-i | \leq \log(N) N^{\eps} } \int_{I} \frac{1}{ | \lambda_i - z | } \, dE \notag\\
&+ \frac{ N^{\eps/5}}{N} \sum_{ n : |n-i | > \log(N) N^{\eps} } \int_{I} \frac{1}{ | \lambda_i - z | } \, dE,
\end{align}
where $i$ is  an index such that $\scg_i \in I$.  
The first sum is bounded by
\beq\label{e:a11}
 \frac{ N^{\eps/5}}{N} \sum_{ n : |n-i | < \log(N) N^{\eps} } \int_{I} \frac{1}{ | \lambda_i - z | } \, dE \leq C \frac{ N^{ \eps  + \eps_1 } }{N} \left(  | 1 + \log ( \eta ) | \right).
\eeq
The second sum can be estimated using rigidity, and we obtain
\beq\label{e:a12}
\frac{ N^{\eps/5}}{N} \sum_{ n : |n-i | > \log(N) N^{\eps_1} } \int_{I} \frac{1}{ | \lambda_i - z | } \, dE \leq C N^{\eps/5}  \log (N) | I | \leq C  \frac{ N^{ \eps + \eps_1}}{N}. 
\eeq
Combining \eqref{e:a11} and \eqref{e:a12} completes the proof.
\end{proof}

From the previous two lemmas we quickly deduce the following. 
\bel \label{lem:intGbd}
Let $I$, $\alpha$, $\eps$, and $C$ be as in the previous lemma.  If rigidity and delocalization hold for $X$ with parameters $\eps$ and $\alpha$, then
\beq
\int_I | G_{a_1 b_1 } ( z) | \cdots | G_{a_k b_k } (z) | \, dE \leq  C \frac{ N^{k\eps+ \eps_1} }{N} \left( 1  + | \log ( \eta ) | \right) \left(\frac{1}{N \eta } + 1  \right)^{k-1}.
\eeq
\eel

\subsubsection{Differentiation}
We have seen that the quantity
\beq \label{eqn:rFE1}
\int_I r ( F_E ) \, dE
\eeq
is a good approximation for the fluctuations of the eigenvalue $\lambda_i$.  We now estimate derivatives of this quantity with respect to matrix elements. We use the shorthand ${\partial_{bc} = \partial_{X_{bc}} }$ for such derivatives. Let us rewrite
\begin{align}
F_E &= \frac{N}{  \pi } \int_{\rr^2} \left( i f_E (e) \chi' ( \sigma) - \sigma f'_E (e) \chi' ( \sigma ) \right) m_N (e + i\sigma ) \, de \, d\sigma \notag\\
&+ \frac{2 i N } { \pi } \int_{ | \sigma | > \eta_2}  \int_\R f'_E (e)  \del_{\sigma} ( \sigma \chi ( \sigma ) )   \Re [ m_N (e + i\sigma ) ] \, de \, d\sigma = A_E + B_E. \notag
\end{align}
We have the following estimate for these two quantities.
\bel\label{l:AB}
Fix $\alpha \in (0, 1)$ and $\eps, \eps_1 \in (0, 1/2)$ with $\eps < \eps_1/2$. Suppose that $I$ is an interval such that $I \subseteq [ \scg_{\alpha N/2}, \scg_{(1-\alpha/2)N}]$ and $|I| \leq N^{\eps_1}/N$.  Suppose rigidity and delocalization hold for the matrix $X$ with parameters $\alpha$ and $\eps$, and $1 \le k \le 5$. Then for $E \in I$ we have
\beq\label{e:intervalbounds0}
|\del_{bc}^k A_E | \leq CN^{k \eps}, \qquad | \del_{bc}^k B_E | \leq CN^{k \eps+k\delta},
\eeq
and for the integral over $I$ we have
\beq\label{e:intervalbounds}
\int_I | \del_{bc}^k A_E | \leq CN^{k \eps}|I|, \qquad \int_I | \del_{bc}^k B_E | \leq C\frac{N^{k \eps+\eps_1 + (k-1) \delta}}{N},
\eeq
where $C=C(\alpha, \eps, \eps_1)>0$ is a constant.
\eel

\begin{proof}We recall the Green's function differentiation formula\footnote{This is a straightforward consequence of the resolvent expansion to first order. See \cite[(2.3)]{lectures}.}
\beq
\frac{\partial G_{ij}(z)}{\partial H_{kl}} = - G_{ik} (z)G_{lj}(z).
\eeq
Using it, we find
\beq
\left| \del_{bc} N m (e + i \sigma) \right|  \le C \sum_a | G_{ab} G_{ca} |   \leq C \sum_a | G_{ab} |^2 + |G_{ca} |^2 \leq \frac{C}{ \sigma} \left( | G_{bb} (z) | + |G_{cc} (z) | \right),
\eeq
where in the last line we used the Ward identity \cite[(3.16)]{lectures}.

From this and \eqref{e:betterlaw} we see that
\beq
| \del_{bc} A_E | \leq C N^{\eps}.
\eeq
Similarly we obtain
\beq
| \del_{bc}^k N m (z) | \leq \frac{ C}{ \sigma} \left( | G_{bb} (z) | + |G_{cc} (z) | + |G_{bc} (z) | \right)^{k},
\eeq
and so 
\beq
| \del_{bc}^k A_E | \leq C N^{k \eps}.
\eeq
for any $1\le k \le 5$.  This proves the first bound in \eqref{e:intervalbounds0}, and the first bound in \eqref{e:intervalbounds} follows immediately.

Now we consider $B_E$.  We have the bound
\beq \label{eqn:delBE}
| \del_{bc}^k B_E | \leq C \int_{E}^{E+ \eta_1}  \int_{ \sigma > \eta_2 } \frac{1}{ \eta_1} \frac{C}{ \sigma} \left( |G_{bb} (z) | + G_{cc} (z) | + G_{bc} (z) | \right)^k \, de \, d\sigma.
\eeq
From this and \eqref{e:betterlaw} we see that
\beq
| \del_{bc}^k B_E | \leq N^{k \eps + k \delta_2 }.
\eeq
From \eqref{eqn:delBE} and Lemma \ref{lem:intGbd} we have,
\beq \label{eqn:intdelBE}
\int_{I} | \del_{bc}^k B_E | \leq \frac{ N^{\eps + \eps_1}}{N}  \int_{ \sigma > \eta_2} \chi ( \sigma ) \frac{1}{ \sigma ( N \sigma )^{k-1} } \, d\sigma \leq  \frac{ N^{2 \eps+ \eps_1 + (k-1) \delta} }{N}.
\eeq
This completes the proof.
\end{proof}

\begin{proof}[Proof of Lemma \ref{lem:ev}]
The smoothed eigenvalues $\tilde \lambda_i$ were constructed in Section \ref{sec:regdef}, and the first estimate in \eqref{eqn:evreg} was already derived as \eqref{e:evreg1}.

To control the derivatives of $\tilde \lambda_i$, it suffices to control the derivatives of \eqref{eqn:rFE1} using Lemma~\ref{l:AB}.  We start with $k=1$, and recall that by Theorem~\ref{thm:rig}, delocalization and rigidity hold with overwhelming probability for the generalized Wigner matrix $H$.  We have
\begin{align}
\left| \del_{bc} \int_I r ( F_E ) \, dE \right| \leq C \int_I | \del_{bc} A_E | + | \del_{bc} B_E | \, dE \leq \frac{ N^{\eps_1 + \eps}}{N}.
\end{align}
Now for $k=2$ we have
\begin{align}
\left| \del^2_{bc} \int_I r (F_E ) \, dE \right| \leq C \int_I | \del_{bc}^2 A_E | + | \del_{bc}^2 B_E | + ( | \del_{bc} B_E | + | \del_{bc} A_E | )^2 \, dE.
\end{align}
The first term is bounded by 
\beq
\int_I | \del_{bc}^2 A_E | \, dE \leq C \frac{ N^{\eps_1+2 \eps}}{N}.
\eeq
For the second term we use \eqref{e:intervalbounds} to bound
\beq
\int_E | \del_{bc}^2 B_E |  \, dE \leq C \frac{ N^{ \eps_1 +2 \eps+ \delta_2 }}{N}.
\eeq 
For the last term we estimate
\beq
\int_I  ( | \del_{bc} B_E | + | \del_{bc} A_E | )^2 \, dE \leq N^{\eps+\delta_2 } \int_I  | \del_{bc} B_E | + | \del_{bc} A_E |  \, dE \leq C \frac{ N^{ \eps_1 +2 \eps+ \delta_2 }}{N}.
\eeq
Hence we have with overwhelming probability,
\beq
\left| \del_{bc} \int_I r (F_E ) \, dE \right| \leq \frac{ N^{\eps+2 \eps_1 + \delta}}{N}
\eeq
for any $\eps >0$.  In general the $k$th derivative of the quantity \eqref{eqn:rFE1} is bounded by
\beq
\left| \del_{bc}^k \int_I r (F_E ) \, dE \right| \leq C \sum_{M_E} \int_I M_E \, dE
\eeq
where each $M_E$ is a monomial in the terms $|\del_{bc}^n B_E |$ and $| \del_{bc}^m A_E|$.  There are finitely many monomials in the sum.  If a monomial has $l \leq k$ terms and if we denote the order of the derivative in the $i$th term in the monomial by $n_i$ then $n_1 + n_2 + \cdots + n_l = k$.  We have
\beq
 \int_I M_E \, dE \leq N^{\eps + (k-n_1) \delta_2 } \int_I | \del^{n_1} X_E | \, dE
\eeq
where $X_E$ is either $A_E$ or $B_E$.  We see that
\beq
\int_I | \del^{n_1} X_E | \, dE \leq \frac{ N^{ \eps + \eps_1 + (n_1-1) \delta_2 }}{N}.
\eeq
Hence,
\beq
 \int_I M_E \, dE \leq C N^{2 \eps + \eps_1 + (k-1) \delta_2}{N}.
\eeq
This proves eigenvalue derivative bounds for $H$. We now extend them to the perturbations $\theta_{ab} H$.

Fix $\kappa>0$ and consider $\eta \ge N^{-1 +\kappa}$. From Theorem~\ref{thm:semi} and a resolvent expansion to high order,\footnote{See \cite[(2.3)]{lectures}.} it is straightforward to prove that
\beq\label{e:conc1}
\sup_{  0 \leq \theta \leq 1 } \left| \frac{1}{N} \tr \frac{1}{H-z} - \frac{1}{N} \tr \frac{1}{ \theta_{ab} H-z} \right| \leq \frac{N^{\eps}}{N \eta}
\eeq
and
\beq\label{e:conc2}
\sup_{ 0 \leq \theta \leq 1 }\left| \left( \frac{1}{H-z} \right)_{ij} - \left( \frac{1}{ \theta_{ab} H -z } \right)_{ij} \right| \leq \frac{N^{\eps}}{\sqrt{N}}
\eeq
for any $\eps >0$, with overwhelming probability.  Hence, by standard arguments,\footnote{See \cite[Theorem 2.10]{lectures} and \cite[Section 9]{lectures}.} rigidity and delocalization hold simultaneously for $\theta_{ab} H$ for all choices of $\theta$, and we can apply the above calculations, which give the second estimate of \eqref{eqn:TF}.

Finally, the almost sure bounds in \eqref{e:tildeas} follow from using the trivial bound $|G(z)|_{ij} \le \eta^{-1}$ in the above calculations.
\end{proof}

\section{Extension to maximal gap over an interval} \label{sec:int}

In this section, we extend the method of the previous subsection to cover statistics of the form
\beq
\max_{\lambda_i \in I} \lambda_{i+1} - \lambda_i ,
\eeq
where $I = [a, b]$ is compact subinterval of $(-2, 2)$, which proves Theorem \ref{thm:hatT}.  The main issue is that the indicator function $\1_{ \{ \lambda_i \in I \} }$ is not differentiable.  We pass to a smoothed object with the help of the following  Wegner-type estimate.
\bel \label{lem:weg}
Let $H$ by a generalized Wigner matrix and $E \in (-2 + \kappa, 2 - \kappa)$ for fixed $\kappa >0$.  For every $\eps >0$ there is $\delta  = \delta(\eps)>0$ such that
\beq
\pp \left[ \exists i : | \lambda_i - E | \leq 2 N^{-1-\eps } \right] \leq N^{- \delta}.
\eeq
\eel
\begin{proof} This is a corollary of fixed energy universality for generalized Wigner ensembles with an effective polynomial rate of convergence.  This was proved in \cite{landon2019fixed}.  Alternatively, this estimate was proved for Gaussian divisible ensembles in Section 7 of \cite{LY}.  The estimate can then be transferred to all ensembles using the four moment approach of \cite{TV10,TV11} or the matrix continuity estimate of \cite{QUE}. \end{proof}

Fix a small $\eps_{w} >0$ with corresponding $\delta_w$ as in the statement of the lemma.  Let $\rho$ be a smooth test function such that $\rho(x) = 1$ for $|x| \leq \frac{1}{2}$ and $\rho (x) = 0$ for $|x| \geq 1$.  Let $r(x) : \rr_{\geq 0} \to \rr_{\geq 0}$ be a smooth function that is $1$ for $x \leq \frac{1}{2}$ and $0$ for $x \geq 1$.  Let $i_0$ and $j_0$ be the indices of the closest classical eigenvalues to $a$ and $b$.  Fix a small $\eps_r >0$ and consider the index set $J_{r} = \llbracket i_0 - N^{\eps_r}, j_0 + N^{\eps_r}\rrbracket$.  Consider the function
\beq
f_1 ( \blam)   = r \left( \sum_{ i \in J_{r}}  \rho ( N^{1+\eps_w} ( \lambda_i  -a ) ) \right),
\eeq
define $f_2 ( \blam )$ similarly  but with $b$ instead of $a$.  Note that on the event that no eigenvalue is within distance $N^{-1-\eps_w}$ of $a$ and $b$ we have that $f_1(x)=f_2 (x)=1$.  Moreover, if there is an eigenvalue with distance $N^{-1-\eps_w}/2$ of $a$, then  $f_1$ is zero, and similarly for $b$ and $f_2$.


Let $\chi (x)$ be a smoothed out step function with $\chi (x) = 1$ for $x \leq - \frac{1}{2}$ and $\chi (x) = 0$ for $x \geq \frac{1}{2}$.  Consider the two functions
\beq
g_1 (x) = \chi ( 10 N^{1 + \eps_w } (a - x ) ),\qquad g_2 (x) = \chi ( 10 N^{1+\eps_w} ( x - b) ).
\eeq
The important observation is that if $|\lambda_i - b| > N^{-1-\eps_w}$, then $g_2 ( \lambda_i ) = \1_{ \{x < b \} }$, and moreover that if $g_2^{(k)} ( \lambda_i ) \neq 0$ for any $k \ge 1$, then $f_2$ and all of its derivatives are $0$.  There is a similar consideration with $g_1$ and $f_1$.

We now fix $\beta = N^{\gamma}$ and consider
\beq
\widehat{F} (\blam ) = f_1 ( \blam) f_2 ( \blam) \frac{1}{\beta} \log \widehat{Z},
\eeq
where
\begin{align}
\widehat{Z}(\blam)  &= \sum_{ i \in J_r } g_1 ( \lambda_i ) g_2 ( \lambda_i ) \exp \left[ \nu \beta ( \lambda_{i+1} - \lambda_i ) \right].
\end{align}

\bel Let $I = [a,b] \subset (-2,2)$. Let $\eps_w$, $\eps_r$, and $\delta_w$ be as above.  Let $f_1$, $f_2$ and $g_2$ and $g_1$ be as above.  Let $H$ be a generalized Wigner matrix with eigenvalues $\lambda_i$.  Let $\tillam_i$ be the regularized eigenvalues from Lemma \ref{lem:ev}, with $\eps$ and $\delta$ given.  Choose these parameters to satisfy
\beq \label{eqn:paramchoice}
\eps_w < \delta - \eps, \qquad \eps < \eps_r /2.
\eeq
   Then with probability at least $1 - CN^{-\delta_w}$, we have the estimate,
\beq
\left|  \max_{ i  : \lambda_i \in I }  \nu ( \lambda_{i+1} - \lambda_i )  - \widehat{F} ( \btlam)  \right| \leq C \nu \frac{N^{\eps}}{N^{1+\delta}}+  \frac{ \log(N)}{N^{\gamma}}.
\eeq
\eel
\begin{proof} By Lemma \ref{lem:weg} there is an event with probability at least $1 - 2N^{-\delta_w}$ on which there is no eigenvalue $\lambda_i$ within distance $2 N^{-1-\eps_w}$ of the interval endpoints $a$ or $b$. Further, rigidity \eqref{e:rigidity} holds for a sufficiently small $\eps < \eps_r/2$ with overwhelming probability, as does the first estimate of \eqref{eqn:evreg}. In the latter, we choose $\alpha$ sufficiently small so that the classical eigenvalues corresponding to indices in $\llbracket \alpha N, (1- \alpha ) N\rrbracket$ contain an interval which contains $I$.  

When these three events hold,  $\lambda_i \in I$ if and only if $\tillam_i \in I$, $i \in J_{r} = \llbracket i_0 - N^{\eps_r}, j_0 + N^{\eps_r}\rrbracket$, and no $\tillam_i$ is within distance $N^{-1-\eps_w}$ of the interval endpoints $a$ or $b$ (recall \eqref{eqn:paramchoice}).  From this discussion we see first that
\beq
\left| \max_{ i : \lambda_i \in I } \nu ( \lambda_{i+1} - \lambda_i ) - \max_{ i \in J_r: \tillam_i \in I , } \nu ( \tillam_{i+1} - \tillam_i )  \right| \leq C \nu \frac{N^{\eps}}{N^{1+\delta}}
\eeq
when these events hold.

Secondly, we see that on the event in question that $f_1 ( \tillam_i ) = f_2 ( \tillam_i ) = 1$,  and for $ i \in J_r$ that $g_1 ( \tillam_i ) g_2 ( \tillam_i )= \1_{ \{ \tillam_i \in I \} }$.  
Hence, similarly to Lemma \ref{l:entropy}, we see that
\beq
\left|  \max_{ i : \tillam_i \in I } \nu ( \tillam_{i+1} - \tillam_i ) -  \widehat{F} ( \btlam ) \right| \leq \frac{ \log(N)}{\beta}.
\eeq
This yields the claim. \end{proof}

We now extend Lemma \ref{l:fpartial} to $\widehat{F}$.  Consider $\widehat{F}$ as a function on $\rr^N$, $\widehat{F} ( \boldsymbol v)$, for $\boldsymbol v \in \rr^N$.  
The key observation is that if $\del^{(k)}_{v_j} (f_1 f_2) \neq 0$ for some $k \ge 1$, then there is no $v_i$ within distance $N^{-1-\eps_w}/2$ of $a$ or $b$.  If the latter holds, then $\del_{v_i} (g_1 g_2) =0$ by construction.  
This allows the estimation of derivatives of $\log ( \widehat{Z})$ without worrying about quantities like $g_1'/g_1$ growing large.  

\bel \label{lem:hatF}
Let $\eps_w$, $\eps_r$, $\widehat{F}$, $\widehat{Z}$ be as above.  Let $\A \subseteq \rr^N$ be defined as
\beq
\A = \{ \boldsymbol v \in \rr^N :  2 |v_i - a | \geq N^{-1-\eps_w} \mbox{ and } 2 |v_i - b | \geq N^{-1 - \eps_w} ,\forall  i \in J_r \}.
\eeq
For every $d \geq 1$, we have
\begin{align}
&\sum_{ j_1, \cdots j_d } | \del_{j_1} \cdots \del_{j_d } ( f_1 (\boldsymbol v) f_2 ( \boldsymbol v) ) | \notag\\
\leq & C_d N^{d \eps_w+d} \left( | \{ i \in J_r : |v_i - a | \leq  N^{-1-\eps_r} \} | + | \{ i  \in J_r: |v_i - b | \leq  N^{-1-\eps_r }\} | \right)^d \1_{\A}
\end{align}
and
\beq
|f_1  f_2 | \leq \1_{ \A}.
\eeq
We have also
\begin{align}
\1_{ \A} \sum_{j_1, \cdots j_d } | \del_{j_1} \cdots \del_{j_d} \log \widehat{Z} | \leq C_d ( \beta \nu)^d.
\end{align}
\eel
\begin{proof} The first two estimates follow from the fact that $f_1 f_2$ is the constant function $0$ on the set $\A^c$, and that $\rho^{(k)} ( x)$ is supported on the set $|x| \leq 1$.  For the proof of the final estimate we observe that by the above discussion that if $\boldsymbol v \in \A$ we have $\del_{v_j}^{(k)} g_1 (v_j) g_2 (v_j) = 0$ and $g_1 (v_j) g_2 (v_j) = \1_{ \{ v_j \in I \}}$.  With this in mind, the calculations in Lemma \ref{l:fpartial} go through without change, and we find the final estimate. \end{proof}

With these preparations, we are ready to prove Theorem \ref{thm:hatT}.

\begin{proof}[Proof of Theorem \ref{thm:hatT}] Rigidity implies that with overwhelming probability, 
\beq
| \{ i \in J_r : |v_i - a | \leq  N^{-1-\eps_r} \} | + | \{ i  \in J_r: |v_i - b | \leq  N^{-1-\eps_r }\} |  \leq N^{\eps}
\eeq
for any $\eps >0$.  From Lemma \ref{lem:hatF} and Lemma \ref{lem:ev} we find the analog of the estimate \eqref{eqn:TF}, 
\beq
\sup_{\theta} | \del_{ij}^d \widehat{F} ( \theta_{ab} H )  | \leq \frac{C_d  N^{d \delta}}{\beta} \left( 1 + N^{d ( \gamma +\eps - \delta )} + N^{2 \eps + d\eps_w } \right) ,
\eeq
which holds with overwhelming probability. We also find the same almost sure estimate \eqref{e:as}.  In the case $\ell=1$, the theorem is proven using these estimates, similarly to how Theorems \ref{thm:4mm} and \ref{thm:ito} are deduced from Proposition \ref{prop:main}.     The extension to $\ell >1$ requires only replacing $\widehat{Z}$ by quantities like
\beq
\sum_{i_1, \cdots i_\ell} \left(  \prod_{j} g_1 ( \tillam_{i_j} ) g_2 ( \tillam_{i_j} ) \right) \exp \left[ \beta \nu \sum_{k=1}^\ell \tillam_{i_k+1} - \tillam_{i_k} \right]
\eeq
and proceeding as before.
\end{proof}

\section{Universality corollaries} \label{sec:univ}

We consider the two processes
\beq
\d x_i  = \sqrt{\frac{2}{N \beta }}\, d B_i + \frac{1}{N} \sum_{j \neq i} \frac{1}{ x_i - x_j}\, \d t
\eeq
and
\beq
\d y_i = \sqrt{ \frac{2}{ N \beta}}\, d B_i + \frac{1}{N} \sum_{j \neq i } \frac{1}{ y_i -y_j}\, \d t.
\eeq
For the initial data we take $x_i (0) = \lambda_i (H)$ for all $i \in \llbracket 1,  N \rrbracket$, where $H$ is a generalized Wigner matrix, and $y_i (0) = \lambda_i (G)$, where $G$ is the Gaussian ensemble of the same symmetry class.  The parameter is $\beta=1$ in the real case and $\beta=2$ in the complex case.  Then the proof of \cite[Theorem 3.1]{landon2019fixed} implies the following result.  We comment on adapting the argument there to the simpler setting here in Appendix \ref{sec:homog}. The lemma  may also be directly cited from \cite[Corollary 3.2]{bourgadeprep}, which appeared after this paper was written.

\bep \label{prop:homo}
Let $t = N^{\omega}/N$, and let $\alpha \in ( 0 ,1/2)$.  Assume that $0 < \omega < 1/2$.  Then for all $i \in \llbracket \alpha N, (1- \alpha)N\rrbracket$ and all $\eps >0$, we have with overwhelming probability
\beq \label{eqn:homog}
|x_{i+1} (t) - x_i (t) - ( y_{i+1} (t) - y_i (t) ) | \leq  \frac{N^{\eps}}{ N^{1+c}}
\eeq
for some constant  $c=c(\alpha, \omega)>0$ (independent of $\eps$).
\eep

\begin{remark}
The above flow is not quite the same as \eqref{d:OU}. However since $t \ll 1$, Proposition~\ref{prop:homo} holds also for the flow \eqref{d:OU}, as explained in \cite[Section 2.3]{HLY15}. We use this fact without comment in what follows.
\end{remark}
\begin{proof}[Proof of Corollary \ref{cor:univ}]
We first consider the statement for $\mathcal T_{\ell, J}$. By the results of \cite[Section 16.2]{EYbook}, there exists a Gaussian divisible ensemble $X_t$ matching $H$ to four moments, in the sense of \eqref{eqn:mmatch}, with $t$ satisfying the hypotheses of Proposition \ref{prop:homo}. By this proposition, the largest gap among $\{ \lambda_{i+1} - \lambda_i  : i \in J \}$ for $X_t$ is the same as the corresponding quantity for the GOE or GUE up to an error of size $o(N^{-1})$. We deduce that the conclusion of Corollary \ref{cor:univ} holds for $X_t$. Then by Theorem~\ref{thm:4mm}, it holds for $H$ too.

For the quantity $\widehat{\mathcal T}_{\ell, I}$ there is a small issue, which is that even if \eqref{eqn:homog} holds, it may not be true that $x_i (t) \in I$ if and only if $y_i (t) \in I$.  However, by Lemma \ref{lem:weg} we know that for every $\eps_w >0$ there is a $\delta_w >0$ such that the event that there is no $y_i (t)$ within distance $N^{-1-\eps_w}$ of the endpoints of $I$ holds with probability at least $1-N^{-\delta_w}$.  If we take $\eps_w$ small enough so that $N^{-1-\eps_w}$ is larger than the error on the right side of \eqref{eqn:homog}, then we see that on this event, $x_i (t) \in I$ if and only if $y_i (t) \in I$.  Therefore, with probability at least $1-  CN^{-\delta_w}$, for some $\delta_w >0$, the statistics $\widehat{\mathcal T}_{\ell, I}$ of $H_t$ match those of the corresponding Gaussian ensemble.  The remainder of the argument is similar to $\mathcal T_{\ell, J}$ and we conclude using Theorem~\ref{thm:hatT}. \end{proof}

\begin{remark}
We are thankful to the referee for pointing out the following simpler proof of Corollary \ref{cor:univ} for $\widehat{\mathcal T}_{\ell, I}$, which has the advantage of avoiding the need to first prove Theorem~\ref{thm:hatT}. If $I=[a,b]\subset ( -2 ,2 )$, then by the rigidity estimate \eqref{e:rigidity} the random set of indices $J(I) = \{ i : \lambda_i \in I \}$ satisfies 
\beq\label{e:possible}
\llbracket  \gamma_a + N^\eps , \gamma_b - N^{\eps} \rrbracket \subset  J \subset \llbracket  \gamma_a - N^\eps , \gamma_b + N^{\eps} \rrbracket
\eeq
with overwhelming probability, where $\gamma_a$ and $\gamma_b$ are the classical eigenvalue locations closest to $a$ and $b$, respectively, and $\eps >0$ is arbitrary. We may clearly neglect the exceptional set where this does not occur, using the trivial bound $|S| \le \|S \|_\infty$. 

Using Corollary \ref{cor:univ} for ${\mathcal T}_J$, and observing that this estimate is uniform in $J$ as in the statement of Theorem \ref{thm:ito}, we obtain
\beq\label{e:smaller?}
\sup_J \left| \E_{H} {S(\nu \T_{\ell,J}(\blam))} - \E_{G} {S(\nu \T_{\ell,J}(\blam))} \right| \leq C \left(\sup_{0\le d \le 5} \| S^{(d)}\|_\infty \right) N^{-c},
\eeq
where the supremum is taken over all $J$ satisfying \eqref{e:possible}. It then suffices to show
\beq\label{e:suffices11}
\sup_J \left| \E_{G} {S(\nu \T_{\ell,J}(\blam))} - \E_{G} {S(\nu \T_{\ell,I}(\blam))} \right| \leq C \left(\sup_{0\le d \le 5} \| S^{(d)}\|_\infty \right) N^{-c},
\eeq
for some $c>0$ (possibly smaller than the constant in \eqref{e:smaller?}). In turn it suffices to show that there exists $c>0$ such that, with probability at least $ 1 - N^{-c}$, the maximal gap over such $J$ does not correspond to eigenvalues with indices in 
\beq\label{e:possible2}\llbracket \gamma_a - N^\eps, \gamma_a + N^\eps \rrbracket \cup \llbracket \gamma_b - N^\eps, \gamma_b + N^\eps \rrbracket.\eeq 
By the treatment of the display before (3.5) in \cite{maxgapGUE}, there exists a constant $c' >0$ such that the maximum gap over $I$, and hence over the indices $\llbracket  \gamma_a + N^\eps , \gamma_b - N^{\eps} \rrbracket$, is at least $c' \sqrt{\log N} /N$ with probability at least $1 - N^{-c'}$. On the other hand, we see using \cite[Lemma 3.2]{maxgapGUE}, \cite[Lemma 3.5]{maxgapGUE}, and a union bound that the maximal gap among the indices in \eqref{e:possible2} is greater than $c' \sqrt{\log N} /N$ with probability at most $N^{-c''}$, for some constant $c'' (c', \eps)>0$, where we have chosen $\eps$ small enough. This proves \eqref{e:suffices11}, and hence the corollary, after bounding $S$ using the trivial bound $|S| \le \| S \|_\infty$ on the exceptional set where the maximal gap corresponds to indices in \eqref{e:possible2}.
\end{remark}

\begin{proof}[Proof of Corollary \ref{cor:herm}]
The first result of Corollary \ref{cor:herm} will follow from applying Corollary \ref{cor:univ} with $S$ a polynomial, and the corresponding convergence of the largest gaps for the GUE given by \cite[Theorem 1.7]{maxgapGUE}. For notational convenience, set $T = \frac{MN}{\sqrt{32 \log N }} \widehat{\T}_{\ell_N, I}$. We recall $H$ is the generalized Wigner matrix under consideration, and let $G$ be the Gaussian ensemble corresponding to the symmetry class of $X$.

Note that Corollary \ref{cor:univ} may be applied nontrivially for any test function $S$ satisfying $\sup_{ 0\leq d \leq 3} \| S^{(d)}(x)\|_\infty \ll N^c$, where $c$ is as in Corollary \ref{cor:univ}. We let $S$ be an $N^{\varepsilon}$-smoothed regularization of the $L^p$ norm, with $S(x) = |x - 1|^p$ when $|x| < 4 N^\varepsilon$ and $\sup_{ 0\leq d \leq 3} \| S^{(d)}(x)\|_\infty \le CN^{3 p \varepsilon}$.

By Theorem \ref{thm:rig}, the event $\mathcal A$ on which $|\lambda_i - \scg_i| < N^{ -1 + \varepsilon/2}$ for all $i$ such that $\gamma_i \in (a-\varepsilon,b+\varepsilon)$ holds with overwhelming probability.  We now claim that, for any $D > 0$, there exists a $C(D, p, \eps) >0$ such that
\begin{align}\label{64}
\E_H[|T - 1 |^p ] &\leq \E_H[ |T  - 1|^p \one_{\mathcal A}] + \E_H[|T  - 1|^p \one_{\mathcal A^c}] \\
\label{65}&\leq \E_G[|T- 1|^p] + \left( \E_H[ S(T)] - \E_G[S(T)] \right) +C  N^{2p - D}.
\end{align}
In the last inequality, we used the fact that $S(T) = | T - 1|^p$ on the event $\mathcal A$ to split the first term in \eqref{64}. For the second term in \eqref{64}, we used Cauchy--Schwarz and the fact that $\mathcal A$ holds with overwhelming probability to write 
\beq
\E_H[|T  - 1|^p \one_{\mathcal A^c}] \le \sqrt{\E_H[|T  - 1|^{2p} ] \E_H[\one_{\mathcal A^c}]} = N^{-D} \sqrt{\E_H[|T  - 1|^{2p} ] },
\eeq
for $N$ large enough, and bounded $\E_H[|T  - 1|^{2p} ]$ as follows. Observe that $T$ is is bounded by $N$ times the sum of the absolute values of the largest and smallest eigenvalues. The latter quantity can be bounded by (a constant times) the Frobenius norm of $H$, $\| H \|_F$, and we quickly see that \beq \E \| H \|_F \le \sqrt{\E \| H \|^2_F} \le C \sqrt{N}\eeq from the moment growth and independence assumptions on the matrix entries made in Definition \ref{d:gwigner}.

Then in \eqref{65}, the first term is $o(1)$ by \cite[Theorem 1.7]{maxgapGUE}, the second term is $CN^{3 p \varepsilon-c}$ by Corollary \ref{cor:univ}, and the third term is $O(N^{-1})$ when $D$ is chosen large enough. Therefore, it converges to zero.

The proof of the second result is similar. It follows from choosing $S$ to be smooth functions bracketing indicator functions, after a shift and rescaling by some powers of $\log(N)$, which do not affect the polynomial error rate.  For instance, take $S(x)$ such that $S(x) = 1$ for $x \in I_1 = [a_1, b_1]$, $S(x) = 0$ for $x\not\in (a_1 - N^{-\varepsilon}, b_1 + N^{-\varepsilon})$, and  \beq \sup_{x \in (a_1 - N^{-\varepsilon}, a_1) \cup (b_1, b_1 + N^{-\varepsilon})}  |S^{(d)}(x)| \leq C N^{d \varepsilon}\eeq for $0 \leq d \leq 3$.
\end{proof}

\appendix

\section{Eigenvalue Coupling} \label{sec:homog}

The purpose of this appendix is to indicate how Proposition \ref{prop:homo} may be deduced from the existing literature.  It is essentially a consequence of \cite[Theorem 3.1]{landon2019fixed}.  However, this theorem was proved in a more general context and so there are some assumptions that may seem out of place.  We fix the parameter $\alpha >0$ and restrict our considerations to indices $i \in \llbracket \alpha N,(1-\alpha)N\rrbracket$. We also fix $\omega \in (0,1/2)$ and set $t = N^{-1 + \omega}$. The constants below depend on these choices, but we omit this in the notation for this section. (The dependencies are written explicitly in Proposition \ref{prop:homo}.)

In \cite{landon2019fixed}, the behavior of $x_i (t)$ was considered for more general initial data $x_i (0)$ than just the generalized Wigner setting.  One of the key differences is that the empirical eigenvalue density was no longer a semicircle, but given by a different density coming from free probability, called the \emph{free convolution}.  In that work it is denoted by $\rho_{\fc, t}$.  The regularity of this density is controlled by the parameter $t$, and so in \cite{landon2019fixed}, the particles were allowed to evolve until a fixed time $t_0$ before the coupling between $x_i$ and $y_i$ was introduced.  This allowed the density $\rho_{\fc, t}$ to have some regularity before being manipulated.

In the present setting, both $x_i (t)$ and $y_i (t)$ are described by the semicircle density, but supported on the interval $\sqrt{1+t} [-2, 2]$.  The semicircle density is smooth in the bulk, so allowing the particles to evolve until $t_0$ is unnecessary in the present setting.  Errors on the right side of the main estimate of Theorem 3.1 of \cite{landon2019fixed} can then be replaced by setting $\omega_0 =1$.\footnote{Alternatively, we could have allowed the particles to evolve up a fixed time $t_0$ in the present setting, before allowing the coupling between $x_i$ and $y_i$, but we would not gain anything by doing so.}

The result of \cite{landon2019fixed} studies the behavior of $x_i (t)$ near a particle index $i_0$.  In the set-up of \cite{landon2019fixed}, it is possible that $i_0$ was very close to $1$ or $N$, for instance $i_0 = o (N)$.  To account for this, the coupling between $x_i$ and $y_i$ was such that $x_{i_0}$ and $y_{N/2}$ shared the same Brownian motion terms; otherwise it could be possible that $y_{i_0}$ did not share the bulk GOE/GUE statistics, as it could have been close to the edge, and the argument of \cite{landon2019fixed} would have failed.  Aside from being macroscopically separated from the edge, the index $N/2$ was chosen for notational convenience.  Moreover, the process $x_i$ was scaled and shifted so that the particle $x_{i_0}$ was close to $y_{N/2}$ (which is close to $0$), and so that the local particle density near this point was the same as the coupled Gaussian ensemble, this being given by $\rho_{\mathrm{sc}} (0)$.  This is reflected in the assumption about $\rho_{\fc, t}$ and the location of the classical eigenvalue of index $i_0$ stated in Theorem 3.1 of \cite{landon2019fixed}.  In the present setting, we are only interested in indices separated from the edge by a macroscopic distance, and the densities of $x_i (t)$ and $y_i (t)$ match, as they are both given by the semicircle.  We can therefore ignore this assumption in our application of Theorem 3.1, as the proof would go through in the present setting.\footnote{An alternative route to the application of Theorem 3.1 is the following.  Instead of applying it directly to $x_i$ and $y_i$, one can construct a third process $z_i (t)$, also coming from a Gaussian ensemble, so that $z_{N/2}(t)$ is coupled to $x_{i_0}$ and $y_{i_0}$ of the original processes.  The processes $x_{i_0}$ and $y_{i_0}$ can then be shifted and rescaled \emph{by the same constants} so that $x_{i_0}$ and $y_{i_0}$ are close to $z_{N/2}$ and the densities match the semicircle there.  Then Theorem 3.1 can be applied  twice, once to $x_i$ and $z_i$, and the separately to $y_i$ and $z_i$.  Taking the differences between these differences, the $z_i$ dependence drops out and after undoing the scalings, one arrives at an appropriate estimate for $x_i$ and $y_i$.}


Finally, \cite[Theorem 3.1]{landon2019fixed} implies that
\beq\label{e:gapdiff}
x_i(t) - y_i (t) =\sum_{|j - i | \leq t N^{1 + \delta }} \xi_t \left( \frac{i-j}{N} \right) ( x_j (0) - y_j (0 ) ) + O(N^{-1-c_1} )
\eeq
for any $\delta>0$ and an appropriate constant $c_1 = c_1(\delta) >0$, with overwhelming probability.  Here, $\xi_t (x)$ is a smooth function obeying estimates as given in Proposition 3.2 of \cite{landon2019fixed}.  We claim that applying \eqref{e:gapdiff} to $i$ and $i+1$ yields
\begin{align}\label{e:ourconclusion}
&|x_i(t) - y_i (t) - ( x_{i+1}(t) - y_{i+1} (t)) | \\ &\label{e:a3}\le 
\left|  \sum_{|j - i | \leq t N^{1 + \delta } } \xi_t \left(\frac{i-j}{N}\right) ( x_j (0) - y_j (0 ) )  -  \sum_{|j - (i+1) | \leq t N^{1 + \delta } } \xi_t \left(\frac{i+1-j}{N}\right) ( x_j (0) - y_j (0 ) )  \right| \\ &\quad \quad \leq  \frac{N^\eps}{N^{1+c_2}}
\end{align}
for some constant $c_2 >0$ with overwhelming probability. This implies the estimate we need for \eqref{eqn:homog}.

The essential point is that taking differences of gaps weighted by the $\xi_t$ kernel produces a cancellation that improves upon the naive bound of $N^{\eps - 1}$ obtained by applying the rigidity estimate \eqref{e:rigidity} to each gap separately. To show this precisely, we write the differences of the $\xi_t$ terms as an integral of $\partial_x \xi_t$ and use the derivative estimate \cite[(3.12)]{landon2019fixed}.\footnote{This idea has appeared previously in \cite[Section 3.6]{che2019universality}.}
We have 
\beq 
\xi_t \left(\frac{i-j}{N}\right) ( x_j (0) - y_j (0 ) )  - \xi_t \left(\frac{i+1-j}{N}\right) ( x_j (0) - y_j (0 ) )  = ( x_j (0) - y_j (0 ) )\int_{(i-j)/N}^{(i-j+1)/N} \partial_x  \xi_t(x)\, dx.
\eeq
 We now use rigidity and the estimate\footnote{Note that \cite[Proposition 3.2]{landon2019fixed} is missing the factor of $N^{-1}$. The statement of Lemma 3.12 later in the work, from which Proposition 3.2 is derived, is correct.}  $\partial_x  \xi_t(x) \le C N^{-1} (x^2 + t^2)^{-1}$ to obtain
 \beq
 \eqref{e:a3} \le \frac{N^\eps}{N} \frac{1}{N} \int_{-tN^{1+\delta} }^{tN^{1+\delta}} \frac{1}{x^2 + t^2}\, dx \le C\frac{ N^\eps}{N^2 t },
 \eeq
as desired. 


\bibliography{maxbib}{}
\bibliographystyle{abbrv}

\end{document}